\documentclass{article}

\pdfoutput=1
\usepackage{amsmath}
\usepackage{amsfonts}
\usepackage{amscd}
\usepackage{epic}
\usepackage{eepic}
\usepackage{graphicx}

\def\BbR{\mathbb{R}}
\def\BbN{\mathbb{N}}

\def\BbZ{\mathbb{Z}}

\def\e{\varepsilon}
\def\a{\alpha}

\def\c{\gamma}
\def\vp{\varphi}

\def\E{\mathcal{E}}
\def\D{\mathcal{D}}
\def\F{\mathcal{F}}

\def\P{\mathcal{P}}
\def\S{\mathcal{S}}

\def\tE{\widetilde{\mathcal{E}}}
\def\Q{\mathcal{Q}}
\def\L{\mathcal{L}}
\def\tL{\widetilde{\mathcal{L}}}
\def\tR{\widetilde{R}}

\def\tV{\widetilde{V}}

\def\sd#1#2{#1\backslash#2}

\usepackage{amsthm}

\theoremstyle{plain}
\newtheorem{Th}{Theorem}[section]
\newtheorem{Lem}[Th]{Lemma}
\newtheorem{Cor}[Th]{Corollary}
\newtheorem{Prop}[Th]{Proposition}
\newtheorem*{Rem}{Remark}

\theoremstyle{definition}
\newtheorem{Def}[Th]{Definition}
\newtheorem*{Not}{Notation}
\newtheorem{Ex}[Th]{Example}
\newtheorem{Assum}[Th]{Assumption}

\theoremstyle{remark}

\def\thm{\begin{Th}}
\def\endthm{\end{Th}}
\def\lemma{\begin{Lem}}
\def\endlemma{\end{Lem}}
\def\cor{\begin{Cor}}
\def\endcor{\end{Cor}}
\def\prop{\begin{Prop}}
\def\endprop{\end{Prop}}
\def\definition{\begin{Def}}
\def\enddefinition{\end{Def}}
\def\remark{\begin{Rem}}
\def\endremark{\end{Rem}}
\def\example{\begin{Ex}}
\def\endexample{\end{Ex}}
\def\demo{\begin{proof}}
\def\enddemo{\end{proof}}

\def\notation{\begin{Not}}
\def\endnotation{\end{Not}}
\def\assumption{\begin{Assum}}
\def\endassumption{\end{Assum}}

\def\ui#1{^{(#1)}}

\def\supp#1{{\rm supp}(#1)}

\def\<{\langle}
\def\>{\rangle}
\def\1{1\!\!\!\!1}
\def\uS{\underline{S}}
\def\oS{\overline{S}}

\def\be{\begin{equation}}
\def\ee{\end{equation}}

\begin{document}

\title{\bf Standard resistance form on the infinite Sierpinski gasket and its trace onto the half-line $[0, \infty)$.}
\author{\bf Jun Kigami \\ Department of Mathematical Sciences \\Ritsumeikan University}

\maketitle

\begin{abstract}
We construct a compatible sequence of resistance forms on the graphs approximating the infinite Sierpinski gasket and show that the resistance form given by the limit of the compatible sequence induces the Dirichlet form associated with the Brownian motion on the infinite Sierpinski gasket. Moreover, using this characterization of the Dirichlet form, we obtain an explicit expression of the jump kernel associated with the trace process of the Brownian motion onto the half-line $[0, \infty)$.
\end{abstract}

\section{Introduction}\label{INT}

The Sierpinski gasket $K$ is the very first example of a non-trivial fractal for which a canonical diffusion process called the Brownian motion has been constructed. See the pioneering works by Goldstein \cite{Go}, Kusuoka\cite{Kus1} and Barlow-Perkins\cite{BP}. More precisely, there are two kinds of Sierpinski gaskets: the Sierpinski gasket, which is compact, denoted by the SG for short and illustrated in Figure~\ref{APP}, and the infinite Sierpinski gasket, which is unbounded, denoted by the ISG for short and illustrated in Figure~\ref{ISGFig}. \par
The Sierpinski gasket $K$ is defined as the invariant set of a triple of contraction mappings $(F_0, F_1, F_2)$ on $\BbR^2$. Namely, $K$ is the unique non-empty compact set satisfying 
\[
K = F_0(K) \cup F_1(K) \cup F_2(K),
\]
where $F_i: \BbR^2 \to \BbR^2$ is defined as
\[
F_i(z) = \frac 12(z - p_i) + p_i
\]
for any $z \in \BbR^2$ and $i = 0, 1, 2$ with $p_0 = (0, 0), p_1 = (0, 1)$ and $p_2 = (\frac 12, \frac {\sqrt 3}2)$. The Brownian motion is constructed as a proper scaling limit of simple random walks on the graphs $\{V_m\}_{m \ge 0}$, which appear in Figure~\ref{APP}. The Dirichlet form associated with the Brownian motion is called the standard resistance form on the Sierpinski gasket. Roughly, the resistance forms are a special class of Dirichlet forms where the conductance between any two points is always finite and positive. See Appendix~\ref{APD} for the exact definition of resistance forms. To be exact, one needs to modify a resistance form to make a Dirichlet form from it. See \cite[Theorem~9.4]{Ki16} for details. The standard resistance form on the SG is defined as the limit of a compatible sequence of resistance forms as is seen in Theorem~\ref{SRF.thm30}. See Section~\ref{SRF} for more details.\par
The infinite Sierpinski gasket $K\ui{\infty}$ is defined as a union of scaled Sierpinski gaskets $2^nK$, i.e. 
\[
K\ui{\infty} = \bigcup_{n \ge 0} \rho^{(n)}(K),
\]
where $\rho^{(n)}(z) = 2^nz$ for $z \in \BbR^2$. On the ISG, as in the case of the SG, the Brownian motion on the ISG is constructed as a scaling limit of random walks on (infinite) graphs approximating the ISG. (See \cite{Go, Kus1, BP}.) Unlike the case of the SG, the associated Dirichlet form has not been given as a limit of a compatible sequence of resistance forms but as a limit of Dirichlet forms on the scaled Sierpinski gaskets. See \cite{BP, FHK}. In fact, there has been no explicit construction of a compatible sequence of resistance forms deriving the Brownian motion before. The first main purpose of this paper is to give an exact description of the compatible sequence and its limit, called the standard resistance form on the ISG associated with the Brownian motion. More precisely, we will construct a compatible sequence of resistance forms in Theorem~\ref{ISG.thm10} and define the standard resistance form on the ISG as its limit. Then we will identify the Dirichlet form derived from the standard resistance form as the one associated with the Brownian motion in Section~\ref{BMI}. \par
\begin{figure}
\centering
\includegraphics[width = 250pt]{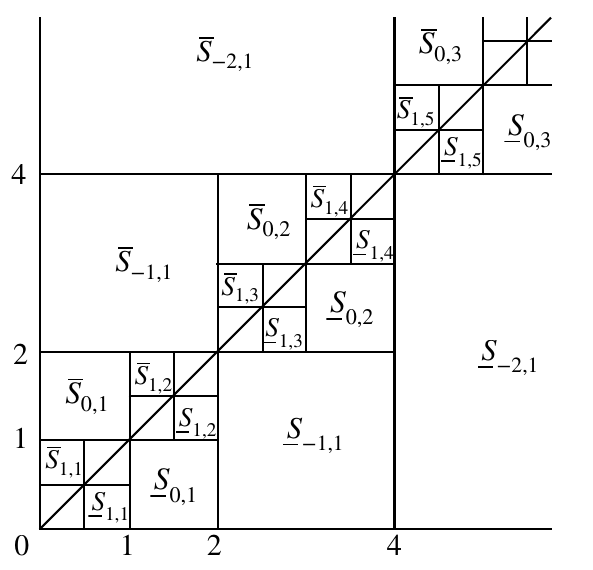}
\caption{Partition of $I\ui{\infty} \times I\ui{\infty} \backslash \{(x, x)| x \in I\ui{\infty}\}$}\label{parI2}
\end{figure}

The second purpose of this paper is to investigate the trace of the Brownian motion on the ISG onto the half line $[0, \infty) \times \{0\}$, which appears as the bottom line of the ISG in Figure~\ref{ISGFig}. For simplicity, we identify $[0, \infty) \times \{0\}$ with $[0, \infty)$ and denote it as $I\ui{\infty}$.
Roughly, the trace process is the process which only observes when the Brownian motion hits $I\ui{\infty}$. In between two consecutive hits, the process stays at the preceding hitting point, so that it is a pure jump process. Our main interest is to obtain an explicit description of the jump kernel associated with the trace. More precisely, define $\nu$ as the restriction of the Lebesgue measure on $\BbR$ to $I\ui{\infty}$, and let $(\E^{(\infty)}_I, \D^{(\infty)}_I)$ be the Dirichlet form on $L^2(I\ui{\infty}, \nu)$ associated with the trace. Then the Dirichlet form has the following expression: 
\[
\E^{(\infty)}_I(u, v) = \int_{I\ui{\infty}}\int_{I\ui{\infty}} (u(x) - u(y))(v(x) - v(y))J^{(\infty)}(x, y)\nu(dx)\nu(dy),
\]
for any $u, v \in \D^{(\infty)}_I$, where $J^{(\infty)}(\cdot, \cdot)$ is called the jump kernel. Before stating our mani result on the jump kernel $J\ui{\infty}$, we need to introduce a partition of $\sd{(I\ui{\infty})^2}{\{(x, x)| x \in I\ui{\infty}\}}$, which is illustrated in Figure~\ref{parI2}.
\prop\label{INT.prop10}
For $i \in \BbN$ and $n \in \BbZ$. Define
\[
\uS_{n, i} = \bigg(\frac{2i - 1}{2^n}, \frac{2i}{2^n}\bigg] \times \bigg[\frac{2i - 2}{2^n}, \frac{2i - 1}{2^n}\bigg),
\]
\[
\oS_{n, i} = \bigg[\frac{2i - 2}{2^n}, \frac{2i - 1}{2^n}\bigg) \times \bigg(\frac{2i - 1}{2^n}, \frac{2i}{2^n}\bigg]
\]
and
\[
S_{n, i} = \uS_{n, i} \cup \oS_{n, i}.
\]
Then $S_{n, i} \cap S_{m, j} = \emptyset$ if and only if $(n, i) \neq (m, j)$ and
\[
\bigsqcup_{i \in \BbN, n \in \BbZ} S_{n, i} = I\ui{\infty} \times I\ui{\infty} \backslash \{(x, x)| x \in I\ui{\infty}\},
\]
where $\bigsqcup$ means a disjoint union.
\endprop

Using this partition, we have the following explicit description of $J\ui{\infty}$.

\thm\label{INT.thm10}
\begin{equation}\label{INT.eq10}
J\ui{\infty}(x, y) = \frac{35}{16}\cdot\frac{14}{17}\Big(\frac {20}3\Big)^{n - 1}
\end{equation}
if $(x, y) \in \bigcup_{i \in \BbN}S_{n, i}$.
\endthm

Such an exact expression of the jump kernel has already been obtained in \cite{KiKTaka} for the case of the SG.\par
The organization of this paper is as follows. In Section~\ref{SRF}, we review the construction of the standard resistance form on the SG. In Section~\ref{ISG}, we construct the standard resistance form on the ISG through a compatible sequence of resistance forms on the graphs approximating the ISG. In Section~\ref{BMI}, we identify the Dirichlet form derived from the standard resistance form on the ISG as the one which is associated with the Brownian motion on the ISG. Section~\ref{TRI} is devoted to studying the Dirichlet form associated with the trace of the Brownian motion onto $I\ui{\infty}$. In particular, we will give a proof of Theorem~\ref{INT.thm10}. Finally, we briefly review the theory of resistance forms.\par
In this paper, to represent a quadratic form, we use the diagonal part of it. More precisely, let $\Q(\cdot, \cdot)$ be a quadratic form, i.e. a bilinear form. Define its diagonal part $\Q(\cdot)$ by $\Q(u) = \Q(u, u)$. By the parallelogram law, it follows that
\[
\Q(u, v) = \frac 14(Q(u + v) - Q(u - v))
\] 
Hence, to know the diagonal part $\Q(\cdot)$ is to know the quadratic form $\Q(\cdot, \cdot)$.

\setcounter{equation}{0}
\section{The standard resistance form on the Sierpinski gasket}\label{SRF}
In this section, we briefly review the construction and basic properties of the standard resistance form $(\E, \F)$ on the Sierpinski gasket $K$. First, we define the Sierpinski gasket. The following theorem is a special case of the fundamental theorem on the existence of self-similar sets. See \cite[Chapter 1]{AOF} for example.

\begin{figure}
\centering
\includegraphics[width = 300pt]{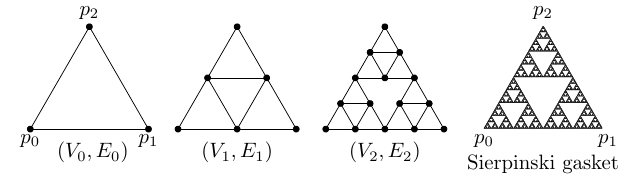}
\caption{Sierpinski gasket and associated graphs}\label{APP}
\end{figure}

\thm\label{SRF.thm10}
Let $p_0$, $p_1$ and $p_2$ be points in $\BbR^2$ given as $p_0 = (0, 0), p_1 = (1, 0)$ and $p_2 = (\frac 12, \frac {\sqrt 3}2)$. Define $F_i:\BbR^2 \to \BbR^2$ for $i = 0, 1, 2$ by 
\[
F_i(z) = \frac 12(z - p_i) + p_i
\]
for $z \in \BbR^2$. Then there exists a unique non-empty compact set $K$ such that
\[
K = F_0(K) \cup F_1(K) \cup F_2(K).
\]
The compact set $K$ is called the Sierpinski gasket, which is illustrated in the rightmost figure of Figure~\ref{APP}.
\endthm

Normally, the Sierpinski gasket is equipped with the restriction of the Euclidean metric on $\BbR^2$, which is denoted by $d_*$. \par
The following is a collection of standard definitions to describe further properties of the Sierpinski gasket.

\definition\label{SRF.def10}
For $m \ge 0$, define
\[
W_m = \{0, 1, 2\}^m = \{w_1\ldots{w_m}| w_j \in \{0, 1, 2\}\,\,\text{for $j = 1, \ldots, m$},\}
\]
and
\[
W_* = \bigcup_{m \ge 0} W_m,
\]
where $W_0 = \{\emptyset\}$. Define
\[
V_0 = \{p_0, p_1, p_2\}.
\]
For $m \ge 0$ and $w = w_1\ldots{w_m} \in W_m$, define
\[
F_w = F_{w_1}{\circ}\ldots{\circ}F_{w_m}, K_w = F_w(K), B_w = F_w(V_0), 
\]
\[
V_m = \bigcup_{w \in W_m} B_w\quad\text{and}\quad V_* = \bigcup_{m \ge 0} V_m,
\]
where $F_{\emptyset}$ is the identity map.
\enddefinition

It is easy to see that $\{V_m\}_{m \ge 0}$ is monotonically increasing, i.e. $V_m \subseteq V_{m + 1}$ for any $m \ge 0$. See Figure~\ref{APP}. Moreover, $V_*$ is dense in $K$ with respect to the metric $d_*$.\par
Next we introduce a canonical measure $\mu_*$ on the Sierpinski gasket $K$.

\thm\label{SRF.thm20}
There exists a Borel regular probability measure $\mu_*$ on $(K, d_*)$ satisfying
\[
\mu_*(K_w) = \Big(\frac 13\Big)^{|w|}
\]
for any $w \in W_*$. Moreover, the measure $\mu_*$ coincides with the normalized $\frac{\log 3}{\log 2}$-dimensional Hausdorff measure of $(K, d_*)$. In particular,
\[
\dim_H(K, d_*) = \frac{\log 3}{\log 2},
\]
where $\dim_H(K, d_*)$ is the Hausdorff dimension of $(K, d_*)$.
\endthm

\remark
The word ``normalised'' means that if $\mathcal{H}(\cdot)$ is the $\frac{\log 3}{\log 2}$-dimensional Hausdorff measure on $(K, d_*)$, then
\[
\mu_*(A) = \frac{\mathcal{H}(A)}{\mathcal{H}(K)}
\]
for any Borel set $A \subseteq K$.
\endremark

Now we start to define the standard resistance form $(\E, \F)$ on the Sierpinski gasket as a limit of a compatible resistance forms $\{(V_m, \E_m)\}_{m \ge 0}$. See Appendix~\ref{APD} for a brief review on resistance forms, including the definitions of a resistance form and a compatible sequence. To begin with, we introduce a natural graph structure to the set $V_m$. Let
\[
E_m = \{(F_w(p_i), F_w(p_j))| i, j \in \{0, 1, 2\}, i \neq j, w \in W_m\}
\]
for $m \ge 0$. Then $(V_m, E_m)$ is a non-directed graph. The following theorem follows from \cite[Proposition~3.1.3]{AOF} and \cite[Example~3.1.5]{AOF}.

\notation
For a set $A$, we define
\[
\ell(A) = \{f| f: A \to \BbR\}
\]
\endnotation

\thm\label{SRF.thm30}
For $m \ge 0$, define a quadratic form $\E_m$ on $\ell(V_m)$ by
\[
\E_m(u) = \frac 12\Big(\frac 53\Big)^m\sum_{(x, y) \in E_m}(u(x) - u(y))^2
\]
for $u \in \ell(V_m)$. Then $(\E_m, \ell(V_m))$ is a resistance form on $V_m$ and the sequence $\L = \{(V_m, \E_m)\}_{m \ge 0}$ is a compatible sequence of resistance forms. 
\endthm

By Theorem~\ref{BRF.thm10}, the limit of $\L$, $(\E_{\L}, \F_{\L})$, is a resistance form on the completion of $V_{\L} = V_*$, which is identified with the Sierpinski gasket $K$ by the results in \cite[Chapter~3]{AOF}, in particular \cite[Example~3.1.3 and Theorem~3.3.4]{AOF}.

\definition\label{SRF.def20}
The resistance form $(\E_{\L}, \F_{\L})$ on $K$ is called the standard resistance form on the Sierpinski gasket, and it is denoted by $(\E, \F)$ for simplicity. Furthermore, the resistance metric associated with $(\E, \F)$ is denoted by $R$.
\enddefinition

The following facts, except \eqref{TSS.eq10}, are also the conclusion of the results in \cite[Chapter~3]{AOF}. See \cite{KiKTaka} for \eqref{TSS.eq10}.

\thm\label{SRF.thm40}
Let $C(K) = \{f| f: K \to \BbR, \text{$f$ is continuous on $K$}\}$. Then
\[
\F = \{f| f \in C(K), \lim_{m \to \infty} \E_m(f|_{V_m}) < +\infty\}
\]
and
\[
\E(f) = \lim_{m \to \infty} \E_m(f|_{V_m})
\]
for any $m \ge 0$. Moreover, there exist $c_1, c_2 > 0$ such that
\begin{equation}\label{TSS.eq10}
c_1|x - y|^{\a} \le R(x, y) \le c_2|x - y|^{\a}
\end{equation}
for any $x, y \in K$, where $\a = \frac{\log 5 - \log 3}{\log 2}$.
\endthm

\setcounter{equation}{0}
\section{The standard resistance form on the infinite Sierpinski gasket}\label{ISG}
In this section, we introduce the infinite Sierpinski gasket and construct a resistance form on it as a natural extension of the standard resistance form $(\E, \F)$ on the Sierpinski gasket $K$ given in the previous section.
 
\begin{figure}
\centering
\includegraphics[width = 300pt]{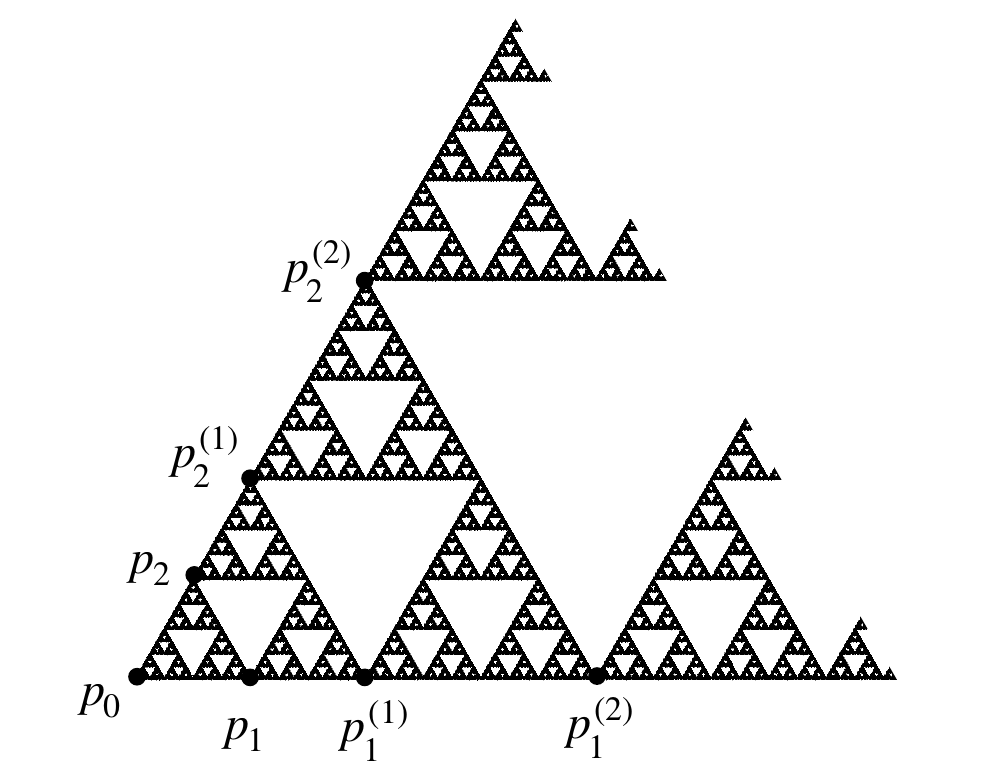}
\caption{The infinite Sierpinski gasket $K^{(\infty)}$}\label{ISGFig}
\end{figure}

 \definition\label{ISF.def10}
For $m \in \BbZ$ and $i \in S$, define $\rho_i^{(m)}: [0, \infty)^2 \to [0, \infty)^2$ by
\[
\rho_i^{(m)}(z) = 2^m(z + p_i)
\]
for $z \in [0, \infty)^2$. For $m \ge 0$ and $n \ge -m$, define
\[
V_n^{(m)} = \rho_0^{(m)}(V_{m + n}), E^{(m)}_n = \{(\rho_0^{(m)}(p), \rho_0^{(m)}(q))| (p, q) \in E_{m + n}\}
\]
and 
\[
\E_n^{(m)}(u) = \Big(\frac 35\Big)^m\E_{m + n}(u{\circ}\rho_0^{(m)})
\]
for $u \in \ell(V^{(m)}_n)$. Moreover, define
\[
p^{(m)}_1 = \rho_0^{(m)}(p_1)\quad\text{and}\quad p^{(m)}_2 = \rho_0^{(m)}(p_2).
\]
In particular, $p_i^{(0)} = p_i$ for $i = 1, 2$. Define 
\[
K^{(m)} = \rho_0^{(m)}(K)\quad\text{and}\quad K\ui{\infty} = \bigcup_{m \ge 0} K\ui{m}.
\]
The unbounded set $K\ui{\infty}$ is called the infinite Sierpinski gasket (ISG for short), which is illustrated in Figure~\ref{ISGFig}. 
\enddefinition

\remark
\[
\E_n^{(m)}(u) = \frac 12\Big(\frac 53\Big)^n\sum_{(p, q) \in E^{(m)}_n} (u(p) - u(q))^2
\]
\endremark

\remark
For a self-similar set like the Sierpinski gasket, we have infinitely many ways to blow it up and obtain its infinite version. See \cite{St0} for example. Though $K\ui{\infty}$ is one of such blowups of the Sierpinski gasket, we call it ``the'' infinite Sierpinski gasket in this paper for simplicity.
\endremark

The following proposition, which is not really needed until the next section, describes basic properties of the infinite Sierpinski gasket as a metric-measure space. 

\prop\label{ISG.prop10}
Let $d_*$ be the restriction of the Euclidean metric of $\BbR^2$ to the infinite Sierpinski gasket $K\ui{\infty}$. Then $\dim_H(K, d_*) = \frac{\log 3}{\log 2}$. Moreover, define $\mu_*\ui{\infty}$ by
\[
\mu_*\ui{\infty}(A) = \frac{\mathcal{H}(A)}{\mathcal{H}(K)}
\]
for a Borel set $A \subseteq K\ui{\infty}$, where $\mathcal{H}$ is the $\frac{\log 3}{\log 2}$-dimensional Hausdorff measure on $K\ui{\infty}$. Then
\[
\mu_*\ui{\infty}(A) = 3^m\mu_*((\rho_0\ui{m})^{-1}(A))
\]
for any $m \ge 0$ and Borel set $A \subseteq K\ui{m}$. In particular, $\mu_*\ui{\infty}(A) = \mu_*(A)$ for any $A \subseteq K$.
\endprop

Hereafter, we use $\mu_*$ to denote $\mu_*\ui{\infty}$ for simplicity if no confusion may occur.\par

Note that $K^{(m)}$ is an enlarged Sierpinski gasket and $K^{(m)} \subseteq K^{(m + 1)}$ for any $m \ge 0$.\par

\begin{figure}[ht]
\centering
\includegraphics[width = 300pt]{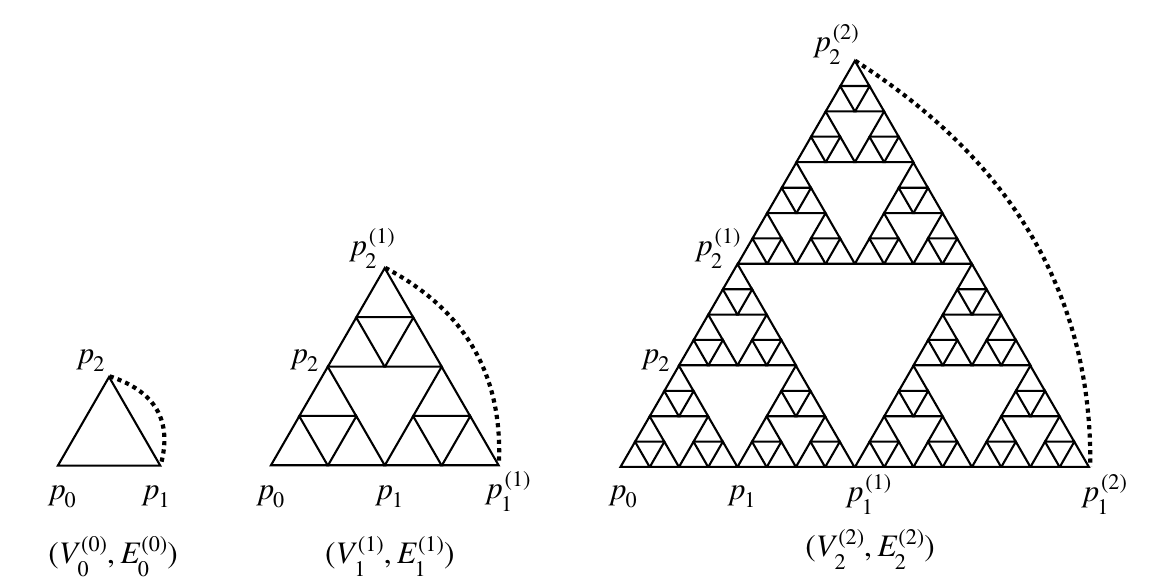}
\vspace{-5pt}
\caption{$(V_m^{(m)}, E_m^{(m)})$ + the edges $(p_1^{(m)}, p_2^{(m)})$ illustrated in dotted lines}\label{Qm}
\end{figure}

\thm\label{ISG.thm10}
For $u \in \ell(V^{(m)}_m)$, define
\[
\Q^{(m)}(u) = \E^{(m)}_m(u) + \frac 56\Big(\frac 35\Big)^m(u(p_1^{(m)}) - u(p_2^{(m)}))^2
\]
Then $\L = \{(V_m^{(m)}, \Q^{(m)})\}_{m \ge 0}$ is a compatible sequence of resistance forms. Moreover, let $(\E^{(\infty)}, \F^{(\infty)}) = (\E_{\L}, \F_{\L})$ and let $R^{(\infty)}(\cdot, \cdot)$ be the resistance metric associated with the resistance form $(\E^{(\infty)}, \F^{(\infty)})$ on $V^{(\infty)} = \bigcup_{m \ge 0} V^{(m)}_m$. Then 
\begin{equation}\label{ISG.eq100}
c_1|x - y|^{\a} \le R^{(\infty)}(x, y) \le c_2|x - y|^{\a}
\end{equation}
for any $x, y \in V^{(\infty)}$, where the constants $c_1$ and $c_2$ are the same as in \eqref{TSS.eq10}. 
\endthm

See Figure~\ref{Qm} for an illustration of the definition of $\Q^{(m)}$, where the edges of $E_m^{(m)}$ are solid lines and the additional edge $(p_1^{(m)}, p_2^{(m)})$ is the dotted line for each $m \in \{0, 1, 2\}$.\par
We will prove this theorem later in this section.\par

By \eqref{ISG.eq100}, the completion of $(V^{(\infty)}, R^{(\infty)})$ is naturally identified with $K\ui{\infty}$. Furthermore, using Theorem~\ref{BRF.thm10}, we may regard $(\E^{(\infty)}, \F^{(\infty)})$ as a resistance form on $K\ui{\infty}$. In what follows, we do so and use $R^{(\infty)}$ to denote the resistance metric associated with the extended $(\E^{(\infty)}, \F^{(\infty)})$. Then the inequality \eqref{ISG.eq100} holds for any $x, y \in K\ui{\infty}$.\par

\definition\label{ISG.def15}
The resistance form $(\E\ui{\infty}, \F\ui{\infty})$ on $K\ui{\infty}$ is called the standard resistance form on the infinite Sierpinski gasket.
\enddefinition

On the enlarged Sierpinski gasket $K\ui{m}$, we have a natural extension of the standard resistance form on the Sierpinski gasket as follows. 

\thm\label{ISG.thm20}
Define
\[
\F^{(m)} = \{u| u: K^{(m)} \to \BbR, u{\circ}\rho_0^{(m)} \in \F\},
\]
and 
\[
\E^{(m)}(u, v) = \Big(\frac 35\Big)^m\E(u{\circ}\rho_0^{(m)}, v{\circ}\rho_0^{(m)})
\]
for $u, v \in \F^{(m)}$. Then $(\E\ui{m}, \F\ui{m})$ is a resistance form on $K\ui{m}$. Let $R\ui{m}$ be the associated resistance metric. Then
\[
c_1|x - y|^{\a} \le R\ui{m}(x, y) \le c_2|x - y|^{\a}
\]
for any $x, y \in K\ui m$, where the constants $c_1$ and $c_2$ are the same as in \eqref{TSS.eq10}. Moreover, if $u \in \F\ui{m + 1}$ and $\supp{u} \subseteq K\ui{m}$, then $u|_{K\ui{m}} \in \F\ui{m}$ and $\E\ui{m + 1}(u) = \E\ui{m}(u|_{K\ui{m}})$.
\endthm

Next, we are going to study the difference between $(\E\ui{m}, \F\ui{m})$ on $K\ui{m}$ and $(\E\ui{\infty}, \F\ui{\infty})$ on $K\ui{\infty}$.

\definition\label{ISG.def20}
For $m \ge n \ge 0$, define
\[
U^{(m, n)}_k = \big(\sd{V_k^{(m)}}{V^{(n)}_k}\big) \cup \{p_1^{(n)}, p_2^{(n)}\}\quad\text{and}\quad E^{(m, n)}_k = \sd{E^{(m)}_k}{E^{(n)}_k}.
\]
Furthermore, for $u \in \ell(U^{(m, n)})$, define
\[
\E^{(m, n)}_k(u ) = \frac 12\Big(\frac 53\Big)^k\sum_{(p, q) \in E^{(m, n)}_k} (u(p) - u(q))^2
\]
and
\[
Q^{(m, n)}_k(u) = \E^{(m, n)}_k(u) + \frac 56\Big(\frac 35\Big)^m(u(p_1^{(m)}) - u(p_2^{(m)}))^2.
\]
\enddefinition

\remark
It follows that $U^{(n, n)}_k = \{p_1^{(n)}, p_2^{(n)}\}$ and
\[
\Q^{(n, n)}_k(u) = \frac 56\Big(\frac 35\Big)^n(u(p_1^{(n)}) - u(p_2^{(n)}))^2.
\]
\endremark

\thm\label{ISG.thm30}
Let$n \ge 0$. The sequence $\L_n = \{(U\ui{m, n}_m, \Q\ui{m, n}_m)\}_{m \ge n}$ is a compatible sequence of resistance forms on finite sets. Let $V\ui{\infty, n} = \bigcup_{m \ge n} U\ui{m, n}_m$ and let $(\E\ui{\infty, n}, \F\ui{\infty, n})$ and $R\ui{\infty, n}$ be the resistance form on $V\ui{\infty, n}$ and the resistance metric associated with the compatible sequence $\L_n$, respectively. Then the completion of $(V\ui{\infty, n}, R\ui{\infty, n})$ is identified with $K\ui{\infty, n} = (\sd{K\ui{\infty}}{K\ui{n}}) \cup \{p_1\ui{n}, p_2\ui{n}\}$. Moreover, 
\[
\F\ui{\infty} = \{u| u: K\ui{\infty} \to \BbR, u|_{K\ui{\infty, n}} \in \F\ui{\infty, n}, u|_{K\ui{n}} \in \F\ui{n}\},
\]
\begin{equation}\label{ISG.eq110}
\E\ui{\infty}(u) = \E\ui{\infty, n}(u|_{K\ui{\infty, n}}) + \E\ui{n}(u|_{K\ui{n}})
\end{equation}
for any $u \in \F\ui{\infty}$,
\begin{multline}\label{ISG.eq120}
\F\ui{\infty} = \{u| u \in C(K\ui{\infty}), \\
\text{$u|_{K\ui{n}} \in \F\ui{n}$ for any $n \ge 0$ and $\lim_{n \to \infty} \E\ui{n}(u|_{K\ui{n}}) < \infty$,}\}
\end{multline}
and
\begin{equation}\label{ISG.eq130}
\E\ui{n}(u|_{K\ui{n}}) \to \E\ui{\infty}(u)
\end{equation}
as $n \to \infty$. 
\endthm

This theorem will be proven later in this section. \par
The final result of this section, the next theorem, shows that compactly supported functions are dense in the domain $\F\ui{\infty}$.

\thm\label{ISG.thm40}
Define
\[
\E\ui{\infty}_0(u) = \E\ui{\infty}(u) + u(p_0)^2
\]
for $u \in \F\ui{\infty}$. Then $\F\ui{\infty} \cap C_0(K\ui{\infty})$ is a dense subset of $\F\ui{\infty}$ with respect to the norm $\sqrt{\E\ui{\infty}_0}$, where $C_0(K\ui{\infty})$ is the collection of continuous functions on $K\ui{\infty}$ having compact supports.
\endthm

The rest of this section is devoted to proofs of Theorems~\ref{ISG.thm10}, \ref{ISG.thm30} and \ref{ISG.thm40}.

\lemma\label{ISG.lemma10}
Let $A_m = \{p_1^{(m)}, p_2^{(m)}\}$. Then 
\[
[Q^{(m + 1, m)}_{m + 1}]_{A_m} = Q^{(m, m)}_{m + 1}
\]
\endlemma

\begin{figure}
\includegraphics[width = 330pt]{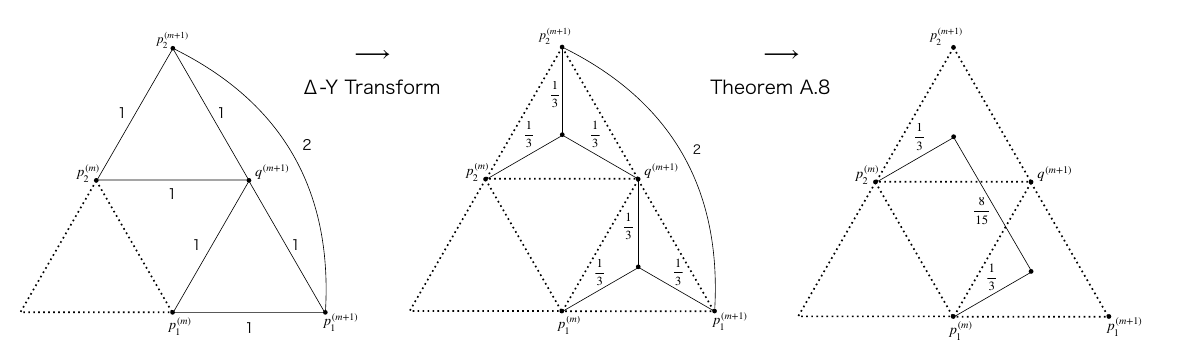}
\vspace{-10pt}
\caption{Modification of $\Big(\displaystyle\frac 53\Big)^m[\Q^{(m + 1, m)}_{m + 1}]_U$}\label{ISG.fig10}
\end{figure}

\demo
Set $q^{(m + 1)} = \frac{p^{(m + 1)}_1 + p^{(m + 1)}_2}2$. Let $U_i = \rho\ui{m}_i(V_{2m + 1})$ and $U_i^0 = \rho\ui{m}_i(V_0)$ for $i = 1, 2$. Then
\[
U^{(m + 1, m)}_{m + 1} = U_1 \cup U_2,
\]
$U_1^0 = \{p^{(m)}_1, p^{(m + 1)}_1, q^{(m + 1)}\}$ and $U_2^0 = \{p^{(m)}_2, q^{(m + 1)}, p^{(m + 1)}_2\}$.
Moreover let $\E_i$ be a resistance from on $U_i$ given by
\[
\tE_i(u) = \Big(\frac 35\Big)^m\E_{2m + 1}(u{\circ}\rho\ui{m}_i).
\]
Then 
\[
\Q_{m + 1}^{(m + 1, m)}(u) = \tE_1(u|_{U_1}) + \tE_2(u|_{U_2}) + \frac 56\Big(\frac 35\Big)^{m + 1}(u(p_1^{(m + 1)}) - u(p_2^{(m + 1)}))^2.
\]
Let $U = U_1^0 \cap U_2^0$. By Lemma~\ref{BRF.lemma10}, 
\begin{multline*}
[\Q_{m + 1}^{(m + 1, m)}]_U(v) \\
= [\tE_1]_{U_1^0}(v|_{U_1^0}) + [\tE_2]_{U_2^0}(v|_{U_2^0})  + \frac 56\Big(\frac 35\Big)^{m + 1}(v(p_1^{(m + 1)}) - v(p_2^{(m + 1)}))^2.
\end{multline*}
Since $[\E_{2m + 1}]_{V_0} = \E_0$ for any $i = 1, 2$, we have
\[
[\Q_{m + 1}^{(m + 1, m)}]_U(v) = \Big(\frac 35\Big)^m\Big(\tE_0^1(u|_{U_1^0}) + \tE_0^2(u|_{U_2^0})
+ \frac 12(v(p_1^{(m + 1)}) - v(p_2^{(m + 1)}))^2\big),
\]
where $\tE_0^i$ is the resistance form on $U_i^0$ given by $\tE_0^i(u) = \E_0(u|_{U_1^0}\circ(\psi_i)^{-1})$ for $i = 1, 2$. Since $U_i^0$ consists of three points, the $\Delta$-Y transform (\cite[Lemma~2.1.15]{AOF}) and Lemma~\ref{BRF.lemma20} yield the desired result. See Figure~\ref{ISG.fig10}, where the figure beside each edge represents the resistance, which is the reciprocal of $C(\cdot, \cdot)$, between corresponding vertices.
\enddemo

\demo[Proof of Theorem~\ref{ISG.thm10}]
First we show $[Q^{(m + 1)}]_{U^{(m)}_m} = Q^{(m)}$. In fact, it follows that
\[
V^{(m + 1)}_{m + 1} = U^{(m + 1, m)}_{m + 1} \cup V^{(m)}_{m + 1}, \quad\quad U^{(m + 1, m)}_{m + 1} \cap V^{(m)}_{m + 1} = A_m
\]
and
\[
Q^{(m + 1)}(u) =  \E^{(m)}_{m + 1}(u|_{V^{(m)}_{m + 1}}) + Q^{(m + 1, m)}_{m + 1}(u|_{U^{(m + 1, m)}_{m + 1}}).
\]
Hence using Lemma~\ref{BRF.lemma10}, we see that
\begin{equation}\label{ISG.eq200}
[Q^{(m + 1)}]_{V^{(m)}_m}(v) = [\E^{(m)}_{m + 1}]_{V^{(m)}_m}(v|_{V_m^{(m)}}) + [Q^{(m + 1, m)}_{m + 1}]_{A_m}(v|_{A_m})
\end{equation}
Since $[\E_{m + 1}]_{V_m} = \E_m$, it follows that $[\E^{(m)}_{m + 1}]_{V^{(m)}_m} = \E^{(m)}_m$. This fact and Lemma~\ref{ISG.lemma10} along with \eqref{ISG.eq200} yield $[Q^{(m + 1)}]_{U^{(m)}_m} = Q^{(m)}$. Consequently, we see that $\L = \{(V^{(m)}_m, Q^{(m)})\}_{m \ge 0}$ is a compatible sequence. Let $x, y \in V^{(m)}_m$. Then $R^{(\infty)}(x, y) = R_{\Q^{(m)}}(x, y)$. Note that
\[
\Q^{(m)}(u) \ge \E^{(m)}_m(u) = \Big(\frac 35\Big)^m\E_{2m}(u{\circ}\rho_0^{(m)}).
\]
Combining this with \eqref{TSS.eq10}, we see that
\begin{multline*}
R_{\Q^{(m)}}(x, y) \le R_{\E^{(m)}_m}(x, y) = \Big(\frac 53\Big)^mR((\rho_0^{(m)})^{-1}(x), (\rho_0^{(m)})^{-1}(y)) \\
\le c_2\Big(\frac 53\Big)^m\Big|\frac x{2^m} - \frac y{2^m}\Big|^{\a} = c_2|x - y|^{\a}.
\end{multline*}
 On the other hand, by \cite[Theorem~A.2]{Ki11}, we have
 \begin{multline*}
R_{\Q\ui{m}}(x, y) = R_{\Q^{(n)}}(x, y)  \ge \bigg(\frac 1{R_{\E^{(n)}_n}(x, y)} + \frac 56\Big(\frac 35\Big)^n\bigg)^{-1}\\
 = \frac{R_{\E^{(n)}_n}(x, y)}{1 + \frac 56\big(\frac 35\big)^nR_{\E^{(n)}_n}(x, y)} \ge \frac{c_1|x - y|^{\a}}{1 + \frac 56\big(\frac 35\big)^nc_2|x - y|^{\a}}.
 \end{multline*}
for any $n \ge m$. Letting $n \to \infty$, we see that $c_1|x - y|^{\a} \le R^{(\infty)}(x, y)$. Thus, we have verified \eqref{ISG.eq100}.
\enddemo

\demo[Proof of Theorem~\ref{ISG.thm30}:Part I]
First we show that $[\Q\ui{m + 1, n}_{m + 1}]_{U\ui{m, n}_m} = \Q\ui{m, n}_m$. Note that
\[
U^{(m + 1, n)}_{m + 1} = U^{(m + 1, m)}_{m + 1} \cup U^{(m, n)}_{m + 1}, \quad\quad U^{(m + 1, m)}_{m + 1} \cap U^{(m, n)}_{m + 1} = A_m
\]
and
\[
Q^{(m + 1, n)}_{m + 1}(u) = Q^{(m + 1, m)}_{m + 1}(u|_{U^{(m + 1, m)}_{m + 1}}) + \E^{(m, n)}_{m + 1}(u|_{U^{(m, n)}_{m + 1}}).
\]
By Lemma~\ref{BRF.lemma10}, it follows that
\begin{equation}\label{ISG.eq300}
[Q\ui{m + 1, n}_{m + 1}]_{U\ui{m, n}_m}(v) = [Q\ui{m + 1, m}_{m + 1}]_{A_m}(v|_{A_m}) + [\E\ui{m, n}_{m + 1}]_{U\ui{m, n}_m}(v|_{U\ui{m, n}_m}).
\end{equation}
On the other hand, we see that
\[
U^{(m, n)}_{m + 1} = \cup_{k = n}^{m - 1} \big(\rho_1^{(k)}(V_{k + m + 1}) \cup \rho_2^{(k)}(V_{k + m + 1})\big)
\]
\[
\E^{(m, n)}_{m + 1}(u) = \sum_{k = n}^{m - 1} \Big(\frac 35\Big)^k\big(\E_{k + m + 1}(u{\circ}\rho_1^{(k)}) + \E_{k + m + 1}(u{\circ}\rho_2^{(k)})\big).
\]
Since $[\E_{k + m + 1}]_{V_{k + m}} = \E_{k + m}$, again using Lemma~\ref{BRF.lemma10}, we obtain
\[
[\E\ui{m, n}_{m + 1}]_{U\ui{m, n}_m}(v) = \sum_{k = n}^{m - 1} \Big(\frac 35\Big)^k\big(\E_{k + m}(v{\circ}\rho_1^{(k)}) + \E_{k + m}(v{\circ}\rho_2^{(k)})\big) = \E\ui{m, n}_m(v).
\]
This along with \eqref{ISG.eq300} and Lemma~\ref{ISG.lemma10} shows $[\Q\ui{m + 1, n}_{m + 1}]_{U\ui{m, n}_m} = \Q\ui{m, n}_m$. Thus we have verified that $\L_n$ is a compatible sequence of resistance forms. \par
Let $x, y \in U^{(m, n)}_m$. Since
\begin{equation}\label{ISG.eq350}
\Q\ui{m}_m(u) = \E\ui{n}_m(u|_{U\ui{n}_m}) + \Q\ui{m, n}_m(u|_{U\ui{m, n}_m})
\end{equation}
and $U\ui{n}_m \cap U\ui{m, n}_m  = A_m$,  Lemma~\ref{BRF.lemma10} implies
\begin{equation}\label{ISG.eq400}
[\Q\ui{m}_m]_{U\ui{m, n}}(v) = [\E\ui{n}_m]_{A_n}(v|_{A_n}) + \Q\ui{m, n}_m(v).
\end{equation}
Using the $\Delta$-Y transform, we deduce
\[
[\E_n]_{\{p_1, p_2\}}(u) = [[\E_n]_{V_0}]_{\{p_1, p_2\}}(u) = [\E_0]_{\{p_1, p_2\}}(u) = \frac 32(u(p_1) - u(p_2))^2.
\]
Hence by \eqref{ISG.eq400}, 
\[
[\Q\ui{m}_m]_{U\ui{m, n}}(v) = \frac 32\Big(\frac 35\Big)^n(u(p_1\ui{n}) - u(p_2\ui{n}))^2 + \Q\ui{m, n}_m(v)
\]
By \cite[Lemma~A.2]{Ki11}, it follows that
\begin{equation}\label{ISG.eq500}
R\ui{\infty}(x, y) \ge \bigg(\frac 1{R\ui{\infty, n}(x, y)} + \frac 23\Big(\frac 53\Big)^n\bigg)^{-1}.
\end{equation}
On the other hand, by \eqref{ISG.eq350},
\begin{multline*}
\inf\{\Q\ui{m}_m(u)|u \in \ell(U\ui{m}_m), u(p_1\ui{n}) = 1, u(p_2\ui{n}) = 0\} \\
\ge \inf\{Q\ui{m, n}_m(u)| u \in \ell(U\ui{m, n}_m), u(p_1\ui{n}) = 1, u(p_2\ui{n}) = 0\}.
\end{multline*}
Hence by \eqref{BRF.eq50}, 
\begin{equation}\label{ISG.eq600}
R\ui{\infty, n}(x, y) \ge R\ui{\infty}(x, y).
\end{equation}
By \eqref{ISG.eq500} and \eqref{ISG.eq600}, a sequence in $V\ui{\infty, n}$ is an $R\ui{\infty}$-Cauchy sequence if and only if it is an $R\ui{\infty, n}$-Cauchy sequence. Therefore, the completions of $V\ui{\infty, n}$ with respect to $R\ui{\infty}$, and $R\ui{\infty, n}$ are the same. By \eqref{ISG.eq100}, the completion is identified with $K\ui{\infty, n}$. \par
Now let $u: V\ui{\infty} \to \BbR$. Then
\begin{equation}\label{ISG.eq610}
\Q\ui{m}_m(u) = \E\ui{n}_m(u|_{U\ui{n}_m}) + \Q\ui{m, n}_m(u|_{U\ui{m, n}_m}).
\end{equation}
Since $\{(U\ui{n}_m, \E\ui{n}_m)\}_{m \ge n}$ and $\{(U\ui{m, n}_m, \Q\ui{m, n}_m)\}_{m \ge n}$ are compatible sequences, it follow that $u \in \F\ui{\infty}$ if and only if $u|_{K\ui{n}} \in \F\ui{n}$ and $u|_{K\ui{\infty, n}} \in \F\ui{\infty, n}$. Letting $m \to \infty$ in \eqref{ISG.eq610}, we obtain
\[
\E\ui{\infty}(u) = \E\ui{n}(u|_{K\ui{n}}) + \E\ui{\infty, n}(u|_{K\ui{\infty, n}}).
\]
Thus we have complete a proof of Theorem~\ref{ISG.thm30} except \eqref{ISG.eq120} and \eqref{ISG.eq130}.
\enddemo

To prove \eqref{ISG.eq120} and \eqref{ISG.eq130}, we need the following lemmas.

\lemma\label{ISG.lemma15}
For any $n \ge 0$, $[\F\ui{\infty}]_{K\ui{n}} = \F^{(n)}$, and
\begin{equation}\label{ISG.eq140}
[\E\ui{\infty}]_{K_n}(v) = \E\ui{n}(v) + \frac 56\Big(\frac 35\Big)^n(v(p_1^{(n)}) - v(p_2^{(n)}))^2
\end{equation}
for any $v \in \F\ui{n}$. Moreover, let $\chi_n: A_n \to \{0, 1\}$ satisfying $\chi_n(p_1\ui{n}) = 1$ and $\chi_n(p_2\ui{n}) = 0$ and let $\psi_n = h_{A_n}^{\E\ui{\infty, n}}(\chi_n)$. Then 
\begin{equation}\label{ISG.eq150}
h_{K_n}^{\E\ui{\infty}}(u)|_{K\ui{\infty, n}} = (u(p_1\ui{n}) - u(p_2\ui{n}))\psi_n + u(p_2\ui{n}).
\end{equation}
\endlemma

\demo
Since $\{(U\ui{m, n}_m, \E\ui{m, n}_m)\}_{m \ge n}$ is a compatible sequence converging the resistance form $(\E\ui{\infty, n}, \F\ui{\infty, n})$, it follows that 
\[
[\E\ui{\infty, n}]_{A_n}(v) = [\Q\ui{m, n}_m]_{A_n}(v) = \Q\ui{n, n}_n(v) = \frac 56\Big(\frac 35\Big)^n(v(p_1\ui{n}) - v(p_2\ui{n}))^2.
\]
On the other hand, for $u \in \F\ui{n}$, we see
\begin{multline*}
[\E\ui{\infty}]_{K_n}(u) = \E\ui{n}(u) + \inf\{\E\ui{\infty, n}(v)| v \in \F\ui{\infty, n}, v|_{A_n} = u|_{A_n}\}\\
= \E\ui{n}(u) + [\E\ui{\infty, n}]_{A_n}(u|_{A_n}).
\end{multline*}
Thus, we have obtained \eqref{ISG.eq140}. Since $h_{A_n}^{\E\ui{\infty, n}}$ is linear, it is straightforward to see \eqref{ISG.eq150}.
\enddemo

\lemma\label{ISG.lemma20}
Let $\tL_m = \{(U\ui{m + 1, m}_k, \E\ui{m + 1, m}_k)\}_{k \ge -m}$. Then $\tL_m$ is a compatible sequence of resistance forms on finite sets. Let $(\E\ui{m + 1, m}, \F\ui{m + 1, m}) = (\E_{\tL_m}, \F_{\tL_m})$. Then $(\E\ui{m + 1, m}, \F\ui{m + 1, m})$ is a resistance form on $K\ui{m + 1, m} = \sd{K\ui{m + 1}}{K\ui{m}} \cup A_m$. Moreover, define 
\begin{multline*}
\tE\ui{m}(u) \\
= \Big(\frac 35\Big)^{m}\Big(\frac 54\sum_{i = 1, 2}(u(p_i\ui{m}) - u(p_i\ui{m + 1}))^2 + \frac 14\sum_{k, l \in \{m, m + 1\}}(u(p_1\ui{k}) - u(p_2\ui{l}))^2\Big)
\end{multline*}
for $u \in \ell(A_m \cup A_{m + 1})$. Then 
\[
[\E\ui{m + 1, m}]_{A_m \cup A_{m + 1}} = [\E\ui{m + 1, m}_k]_{A_m \cup A_{m + 1}} = \tE\ui{m}.
\]
\endlemma

\demo
Define $U_i^k = \rho_i\ui{m}(V_{m + k})$ and $E_i^k = \{(\rho_i\ui{m}(p), \rho_i\ui{m}(q))| (p, q) \in E_{m + k}\}$. for $i = 1, 2$. Then
\[
U\ui{m + 1, m}_k = U_1^k \cup U_2^k,\,\,U_1^k \cap U_2^k = \{q_{m + 1}\}
\]
and
\begin{align*}
\E\ui{m + 1, m}_k(u) &= \Big(\frac 35\Big)^{m}\big(\E_{m + k}(u{\circ}\rho_1\ui{m}) + \E_{m + k}(u{\circ}\rho_2\ui{m})\big)\\
&= \frac 12\Big(\frac 53\Big)^k\big(\sum_{(p, q) \in E_1^k} (u(p) - u(q))^2 + \sum_{(p, q) \in E_2^k} (u(p) - u(q))^2\big)
\end{align*}
Since $[\E_{m + k}]_{V_{m + l}} = \E_{m + l}$ for $-m \le l \le k$, Lemma~\ref{BRF.lemma10} implies
\[
[\E\ui{m + 1, m}_k]_{U\ui{m + 1, m}_l} = \E\ui{m + 1, m}_l.
\]
This shows that $\tL_m$ is a compatible sequence. Let $\tR_m$ be the associated resistance metric on $K$. Let $x, y \in U\ui{m + 1, m}_k$. Then Since $E_1^k \cap E_2^k = \emptyset$, it follows that
\[
\tR_m(x, y) = \begin{cases}
\big(\frac 53\big)^mR((\rho_i\ui{m})^{-1}(x), (\rho_i\ui{m})^{-1}(y))\,&\text{if $x, y \in U_i^k$ for some $i = 1, 2$,}\\
\tR_m(x, q_{m + 1}) + \tR(q_{m + 1}, y)\,&\text{otherwise}.
\end{cases}
\]
In the first case, using \eqref{TSS.eq10} directly, we see that
\begin{equation}\label{ISG.eq900}
c_1|x - y|^{\a} \le \tR_m(x, y) \le c_2|x - y|^{\a}.
\end{equation}
In the second case, using the scaling argument and \eqref{TSS.eq10}, we see that
\begin{equation}\label{ISG.eq910}
c_1(|x - q_{m + 1}|^{\a} + |y - q_{m + 1}|^{\a}) \le \tR_m(x, y) \le c_2(|x - q_{m + 1}|^{\a} + |y - q_{m + 1}|^{\a}).
\end{equation}
On the other hand, an elementary calculation shows that
\[
\frac{\sqrt{3}}2|x - q_{m + 1}| \le |x - y|\quad\text{and}\quad \frac{\sqrt{3}}2|y - q_{m + 1}| \le |x - y|,
\]
Moreover, either $|x - q_{m + 1}| \ge \frac 12|x - y|$ or $|y - q_{m + 1}| \ge \frac 12|x - y|$. By those facts, \eqref{ISG.eq900} and \eqref{ISG.eq910}, we deduce 
\[
c_3|x - y|^{\a} \le \tR_m(x, y) \le c_4|x - y|^{\a},
\]
where $c_3$ and $c_4$ are independent of $x, y$ and $m$. This shows that the completion of $\bigcup_{k \ge -m} U\ui{m + 1, m}_k$ with respect to $\tR_m$ coincides with $K\ui{m + 1, m}$. Finally, let $A = A_m \cup A_{m + 1} \cup \{q_{m + 1}\}$. The $A = U\ui{m + 1, m}_{-m}$ and 
\begin{multline*}
[\E\ui{m + 1, m}]_{A}(u) = [\E^{m + 1, m}_k]_A(u) \\= \E^{m + 1, m}_{-m}(u) = \Big(\frac 35\Big)^{m}\big(\E_{0}(u{\circ}\rho_1\ui{m}) + \E_{0}(u{\circ}\rho_2\ui{m})\big).
\end{multline*}
Now, an elementary calculation shows that the trace of the resistance form in the above equality on $A_m \cup A_{m + 1}$ is $\tE\ui{m}$.
\enddemo

\lemma\label{ISG.lemma30}
Let $\tV\ui{m} = \{p_0\} \cup \bigcup_{k = 0}^{m} A_k$. Then
\begin{equation}\label{ISG.eq700}
[\E\ui{m}_m]_{\tV_m}(v) = \E_0(v|_{V_0}) + \sum_{k = 0}^{m - 1} \tE\ui{k}(v|_{A_k \cup A_{k + 1}}).
\end{equation}
for any $v \in \ell(\tV\ui{m})$. Moreover, 
\begin{multline*}
\widetilde{\F}\ui{\infty} = \{u| u \in C(K\ui{\infty}), \\u|_{K\ui{m}} \in \F\ui{m}\,\,\text{for any $m \ge 0$ and $\lim_{m \to \infty}\E\ui{m}(u|_{K_m}) < \infty$}.\}
\end{multline*}
Then for any $u \in \widetilde{\F}\ui{\infty}$,
\begin{equation}\label{ISG.eq800}
\sum_{m = 0}^{\infty} \Big(\frac 35\Big)^{m}(u(p_1\ui{m})) - u(p_2\ui{m}))^2 < \infty.
\end{equation}
\endlemma

\remark
By \eqref{ISG.eq110}, we have $\F\ui{\infty} \subseteq \widetilde{\F}\ui{\infty}$.
\endremark

\demo
Note that
\[
\E\ui{m}_m(u) = \E_m(u|_{V_m}) + \sum_{i = 0}^{m - 1} \E\ui{i + 1, i}_m(u|_{U\ui{i + 1, i}_m}).
\]
Applying Lemma~\ref{BRF.lemma10} to the above equality and using Lemma~\ref{ISG.lemma20}, we obtain \eqref{ISG.eq700}. Let $u \in \F\ui{\infty}$. Then by \eqref{ISG.eq700}, for any $k \ge 0$, we see that
\[
\lim_{m \to \infty}\E\ui{m}(u|_{K_m}) \ge \E\ui{m}_m(u|_{V\ui{m}_m}) \ge \sum_{k = 0}^{m - 1} \tE\ui{k}(u|_{A_k \cup A_{k + 1}}).
\]
Letting $m \to \infty$, we have \eqref{ISG.eq800} by the definition of $\tE\ui{m}$.
\enddemo

\demo[Proof of Theorem~\ref{ISG.thm30}: Part II]
Now we complete our proof of Theorem~\ref{ISG.thm30}. Namely, we are going to show \eqref{ISG.eq120} and \eqref{ISG.eq130}. Let $u \in \widetilde{\F}\ui{\infty}$. Then for any $m \ge 0$,
\[
\Q\ui{m}(u|_{V\ui{m}_m}) \le \E\ui{m}(u|_{K\ui{m}}) + \frac 56\Big(\frac 35\Big)^m(u(p_1^{(m)}) - u(p_2^{(m)}))^2.
\]
Letting $m \to \infty$ and using \eqref{ISG.eq800}, we see that
\[
\lim_{m \to \infty} \Q\ui{m}(u|_{V\ui{m}_m}) \le \lim_{m \to \infty} \E\ui{m}(u|_{K\ui{m}}) < \infty.
\]
Thus $u \in \F\ui{\infty}$ and $\E\ui{\infty}(u) \le \lim_{m \to \infty} \E\ui{m}(u|_{K\ui{m}})$. On the other hand, by \eqref{ISG.eq110}, it follows that $\E\ui{\infty}(u) \ge \lim_{m \to \infty} \E\ui{m}(u|_{K\ui{m}})$. Thus we have verified \eqref{ISG.eq120} and \eqref{ISG.eq130}. 
\enddemo

\demo[Proof of Theorem~\ref{ISG.thm40}]
Let $u \in \F\ui{\infty}$. Fix $\e > 0$. Then by Theorem~\ref{ISG.thm30} and \eqref{ISG.eq800}, there exists $m$ such that
\[
\E\ui{m, \infty}(u|_{K\ui{\infty, m}}) < \e\quad\text{and}\quad\Big(\frac 35\Big)^{m}(u(p_1\ui{m}) - u(p_2\ui{m}))^2 < \e.
\]
Choose $n > m$ such that
\[
3\Big(\frac 35\Big)^{n}\Big(\frac{u(p_1\ui{m}) + u(p_2\ui{m})}2\Big)^2 < \e.
\]
Let $a_i = u(p_i\ui{m})$ for $i = 1, 2$ and define $v_m:A_{m} \cup A_{m + 1} \to \BbR$ by $v_m(p_i\ui{m}) = a_i$ and $v_m(p_i\ui{m + 1}) = \frac{a_1 + a_2}2$ for $i = 1, 2$ and $v_n: A_n \cup A_{n + 1} \to \BbR$ by $v_n(p_i\ui{n}) = \frac{a_1 + a_2}2$ and $v_n(p_i\ui{n + 1}) = 0$ for $i = 1, 2$.
The define $u_* \in \F\ui{\infty}$ by
\begin{multline*}
u_*|_{K\ui{m}} = u|_{K\ui{m}}, u_*|_{K\ui{m + 1, m}} = h_{A_m \cup A_{m + 1}}^{\E\ui{m + 1, m}}(v_m),\\
u_*|_{K\ui{n, m + 1}} = \frac{a_1 + a_2}2, u_*|_{K\ui{n + 1, n}} = h_{A_n \cup A_{n + 1}}^{\E\ui{n + 1, n}}(v_n),
u_*|_{K\ui{\infty, n + 1}} = 0.
\end{multline*}
Since $\E\ui{k + 1, k}(u_*|_{K\ui{k + 1, k}}) = 0 $ if $k \ge m + 1$ and $k \neq n$, using Lemma~\ref{ISG.lemma20}, we obtain
\begin{multline*}
\E\ui{\infty, m}(u_*|_{K\ui{\infty, n}}) = \E\ui{m + 1, m}(h_{A_m \cup A_{m + 1}}^{\E\ui{m + 1, m}}(v_m)) + \E\ui{n + 1, n}(h_{A_n \cup A_{n + 1}}^{\E\ui{n + 1, n}}(v_n))\\
= \tE\ui{m}(v_m) + \tE\ui{n}(v_n)\\
 = \Big(\frac 35\Big)^{m}(u(p_1\ui{m}) - u(p_2\ui{m}))^2 + 3\Big(\frac 35\Big)^{n}\Big(\frac{u(p_1\ui{m}) + u(p_2\ui{m})}2\Big)^2 < 2e.
\end{multline*}
Therefore, 
\begin{multline*}
\E\ui{\infty}(u - u_*) = \E\ui{\infty, m}((u - u_*)|_{K\ui{\infty, m}}) \\
\le 2\E\ui{\infty, m}(u|_{K\ui{\infty, m}}) + 2\E\ui{\infty, m}(u_*|_{K\ui{\infty, m}}) < 2\e + 4\e = 6\e.
\end{multline*}
Since $u_*(p_0) = u(p_0)$, we have verified the theorem.
\enddemo

\setcounter{equation}{0}
\section{Brownian motion on the infinite Sierpinski \\gasket}\label{BMI}

In this section, we study the Dirichlet form derived from the resistance form $(\E\ui{\infty}, \F\ui{\infty})$. In particular, the primary purpose is to identify the associated diffusion process with the Brownian motion constructed by Barlow-Perkins in \cite{BP} as in the next theorem.\par
 
 \thm\label{BMI.thm10}
 Define
 \[
 \E_1\ui{\infty}(u) = \E\ui{\infty}(u) + \int_{K\ui{\infty}} u^2\mu_*(dz)
 \]
 for $u \in \F\ui{\infty} \cap L^2(K\ui{\infty}, \mu_*)$. Let $\D$ be the closure of $\F\ui{\infty} \cap C_0(K\ui{\infty})$ with respect to the norm $\sqrt{\E_1\ui{\infty}(\cdot)}$. Then $(\E\ui{\infty}, \D)$ is a local regular Dirichlet form on $L^2(K\ui{\infty}, \mu_*)$. Moreover, the diffusion process associated with the local regular Dirichlet form $(\E\ui{\infty}, \D)$ coincides with the Brownian motion constructed by Barlow and Perkins in \cite{BP}. 
 \endthm
 
 \remark
 Originally, the Brownian motion has been constructed on $K\ui{\infty} \cup \pi(K\ui{\infty})$ in \cite{BP}, where $\pi$ is the reflection in the $y$-axis. Identifying $z$ and $\pi(z)$ for $z \in \BbR^2$, we have the counterparts of the results in \cite{BP}, in particular \cite[Theorem~8.1]{BP}, for our case where the space under consideration is $K\ui{\infty}$.
 \endremark
 
 The above theorem enables us to apply the results in \cite{BP} on the properties of the Brownian motion to the diffusion process associated with the Dirichlet form $(\E\ui{\infty}. \D)$. In particular, if $p(t, x, y)$ is the corresponding heat kernel, then it enjoys the following sub-Gaussian heat kernel estimate: there exist $\c_1, \c_2, \c_3, \c_4 > 0$ such that
 \begin{multline}\label{BMI.eq05}
\c_1t^{-\frac{d_S}2}\exp{\Bigg(-\c_2\bigg(\frac{d(x, y)^{d_w}}{t}\bigg)^{1/(d_w - 1)}\Bigg)} \le \\
p(t, x, y) \le \c_3t^{-\frac{d_S}2}\exp{\Bigg(-\c_4\bigg(\frac{d(x, y)^{d_w}}{t}\bigg)^{1/(d_w - 1)}\Bigg)}
\end{multline}
for any $(t, x, y) \in (0, \infty) \times K\ui{\infty} \times K\ui{\infty}$, where the constants $d_S = \frac{\log 9}{\log 5}$ and $d_w = \frac{\log 5}{\log 2}$ are called the spectral dimension and the walk dimension, respectively. On the other hand, this estimate \eqref{BMI.eq05} can also be deduced directly by \cite[Theorem~15.10]{Ki16} without identification with the Brownian motion.\par
 
 The first part of the above theorem is straightforward from \cite[Theorem~9.4]{Ki16}. The rest of this section is devoted to a proof of the second part. The main tool of the proof is the uniqueness of the Brownian motion shown in \cite[Theorem~8.1]{BP}, where restrictions of $(\E\ui{\infty}, \D)$ to unions of neighbouring cells play an important role. To begin with, we define the collection of scaled copies of the Sierpinski gasket in the infinite Sierpinski gasket. 

\definition\label{BMI.def10}
For any $n \in \BbZ$, define
\[
\S_n = \{\rho_0\ui{m}(K_w)| m \ge n, w \in W_{m - n}\}.
\]
For $n \in \BbZ$ and $S \in \S_n$, define $\vp_S$ be the unique map from $\BbR^2$ to $\BbR^2$ satisfying $\vp_S(K) = S$ and $\vp_s(z) = 2^nz + c_S$ for any $z \in \BbR^2$, where $c_S \in \BbR^2$ is independent of $z$.
For $i = 0, 1, 2$, define
\[
p_i(S) = \vp_s(p_i)\quad\text{and}\quad \partial{S} = \{p_0(S), p_1(S), p_2(S)\}.
\]
\enddefinition

It is easy to see that if $S_1 \cap S_2 \neq \emptyset$ for some $S_1, S_2 \in \S_n$, then either $S_1 = S_2$ or $S_1 \cap S_2$ consists of a single point and coincides with $\partial{S_1} \cap \partial{S_2}$.

\definition\label{BMI.def20}
{\rm (1)}\,\,A pair $\P = (S_1, S_2)$ is called a pair of neighbouring cells if $(S_1, S_2) \in \S_n \times \S_n$ for some $n \in \BbZ$ and $S_1 \cap S_2$ consists of a single point. For a pair of neighbouring cells $\P = (S_1, S_2)$, define 
\[
K(\P) = S_1 \cup S_2,\,\,\partial{K(\P)} = \sd{(\partial{S_1} \cup \partial{S_2})}{(\partial{S_1} \cap \partial{S_2})}
\]
and 
\[
K(\P)^o = \sd{K_{\P}}{\partial{K(\P)}}.
\]
{\rm (2)}\,\,For $m, n \in \BbZ$, let $\P = (S_1, S_2) \in \S_n \times \S_n$ and $\P' = (S_1', S_2') \in \S_m \times \S_m$ be pairs of neighboring cells. A bijective map $\vp: K(\P) \to K(\P')$ is called a similitude between $\P$ and $\P'$ if there exist $U_1, U_2 \in O(2)$ such that
\[
\vp(z) = \begin{cases}
2^{m - n}U_1(z - z_0) + z_0'\quad&\text{if $z \in S_1$},\\
2^{m - n}U_2(z - z_0) + z_0'\quad&\text{if $z \in S_2$},
\end{cases}
\]
where $S_1 \cap S_2 = \{z_0\}$ and $S_1' \cap S_2' = \{z_0'\}$.
\enddefinition

\prop\label{BMI.prop10}
Let $\P = (S_1, S_2)$ be a pair of neighboring cells. Define
\begin{multline*}
\D_{\P} = \{u| u: K(\P)^o \to \BbR, \\
\text{there exists $\tilde{u} \in \D$ such that $u = \tilde{u}|_{K(\P)^o}$ and $u|_{(K(\P)^o)^c} \equiv 0$}\}.
\end{multline*}
and
\[
\E_{\P}(u) = \E\ui{\infty}(\tilde{u})
\]
for any $u \in \D_{\P}$, where the function $\tilde{u}$ appears in the definition of $\D_{\P}$. Then $(\E_{\P}, \D_{\P})$ is a local regular Dirichlet form on $L^2(K(\P)^o, \mu_*|_{K(\P)^o})$.
\endprop

Let $X$ be the diffusion process associated with the local regular Dirichlet form $(\E\ui{\infty}, \D)$. Then the diffusion process associated with $(\E_{\P}, \D_{\P})$ is the part process $X^{K(\P)^o}$ killed upon leaving $K(\P)^o$.\par

\thm\label{BMI.thm20}
Let $m, n \in \BbZ$, let $\P \in \S_n \times \S_n$ and $\P' \in \S_m \times \S_m$ be pairs of neighboring cells and let $\vp$ be a similitude between $\P$ and $\P'$. Define $\vp^*(u) = u{\circ}\vp$ for $u \in \D_{\P'}$. Then $\D_{\P} = \vp^*(\D_{\P'})$ and
\[
\E_{\P}(\vp^*(u)) = \Big(\frac 35\Big)^{n - m}\E_{\P'}(u)
\]
for any $u \in \D_{\P'}$.
\endthm

To prove the above theorem, we give an alternative characterization of $(\E_{\P}, \D_{\P})$.

\definition\label{BMI.def30}
Let $n \in \BbZ$ and let $\P \in \S_n \times \S_n$ be a pair of neighboring cells. Choose $l \in \BbN$ such that $K(\P) \subseteq K\ui{l}$. For $k \ge 0$, define
\[
V(\P)_k = V\ui{l}_{k - n} \cap K(\P), E(\P)_k = \{(p, q)| (p, q) \in E\ui{l}_{k - n}, p, q \in V(\P)_k\}
\]
and
\[
\E_{\P, k}(u) = \frac 12\Big(\frac 53\Big)^{k - n}\sum_{(p, q) \in E(\P)_k} (u(p) - u(q))^2.
\]
for $u \in \ell(V(\P)_k)$.
\enddefinition

\remark
The graph $(V(\P)_k, E(\P)_k)$ is independent of the choice of $l$.
\endremark

\lemma\label{BMI.lemma10}
Let $n \in \BbZ$ and let $\P \in \S_n \times \S_n$ be a pair of neighboring cells. Then the sequence $\L(\P) = \{(V(\P)_k, \E_{\P, k})\}_{k \ge 0}$ is a compatible sequence of resistance forms and its limit $(\E_{\L(\P)}, \F_{\L(\P)})$ is a resistance form on $K(\P)$. Moreover, 
\[
\D_{\P} = \{u|_{K(\P)^o} | u \in \F_{\L(\P)}, u|_{\partial{K(\P)}} \equiv 0\}
\]
and $\E_{\P}(u|_{K(\P)^o}) = \E_{\L(\P)}(u)$ for any $u \in \F_{L(\P)}$.
\endlemma

\demo
Let $\P = (S_1, S_2)$. For $i = 1, 2$, define $V(\P)_k^i = V(\P)_k \cap S_i$,
\[
E(\P)_k^i = \{(p, q)| (p, q) \in E(\P)_k, p, q \in S_i\}
\]
and
\[
\E_{\P, k}^i(u) = \frac 12\Big(\frac 53\Big)^{k - n}\sum_{(p, q) \in E(\P)_k^i} (u(p) - u(q))^2.
\]
for $u \in V(\P)_k^i$. Then
\[
V(\P)_k = V(\P)_k^1 \cup V(\P)_k^2, E(\P)_k = E(\P)_k^1 \cup E(\P)_k^2
\]
and
\[
\E_{\P, k}(u) = \E_{\P, k}^1(u|_{V(\P)_k^1}) + \E_{\P, k}^2(u|_{V(\P)_k^2}).
\]
Since $\{(V(\P)_k^i, \E_{\P, k}^i)\}_{k \ge 0}$ is a compatible sequence for each $i = 1, 2$, so is $\{(V(\P)_k, \E_{\P, k})\}_{k \ge 0})\}$. Now choose $l$ such that $K(\P) \subseteq K\ui{l}$. Let $u \in \D_{\P}$ and choose $\tilde{u} \in \D$ as in the definition of $\D_{\P}$. Using Theorem~\ref{ISG.thm30} and the definitions of $\E\ui{l}$ and $\E_{\P, k}$, we see that
\begin{multline*}
\E_{\P}(u) = \E\ui{\infty}(\tilde{u}) = \E\ui{l}(\tilde{u}|_{K\ui{l}}) = \lim_{k \to \infty} \E\ui{l}_{k - n}(\tilde{u}|_{V\ui{l}_{k - n}}) \\= \lim_{k \to \infty} \E_{\P, k}(u|_{V(\P)_k}) = \E_{\L(\P)}(u).
\end{multline*}
Reversing the order of arguments, we obtain the converse direction.
\enddemo

\demo[Proof of Theorem~\ref{BMI.thm20}]
Since $\vp$ maps $V(\P)_k$ and $E(\P)_k$ bijectively to $V(\P')_k$ and $E(\P')_k$, respectively, it follows that
\[
\E_{\P, k}(u) = \Big(\frac 35\Big)^{n - m}\E_{\P', k}(\vp^*(u)).
\]
Letting $k \to \infty$, we obtain the desired statements by Lemma~\ref{BMI.lemma10}.
\enddemo

\demo[Proof of Theorem~\ref{BMI.thm10}]
What we should do is to verify the hypotheses $\rm (H_1), (H_2)$ and $\rm (H_3)$ in \cite[Section~8]{BP} for the diffusion process associated with the Dirichlet form $(\E\ui{\infty}. \D)$. Then the theorem follows from \cite[Theorem~8.1]{BP}. Let $\P = (S_1, S_2) \in \S_n \times \S_n$ be a pair of neighboring cells. In \cite{BP}, the outer triangle of $S_i$ is denoted by $\Delta_i$ for $i = 1, 2$. Moreover, using our terminology, we see that
\[
\Pi(\Delta_1 \cup \Delta_2) = \{\vp| \text{$\vp$ is a similitude between $(S_1, S_2)$ and $(S_1, S_2)$}\},
\]
where $\Pi(\Delta_1 \cup \Delta_2)$ is defined in \cite{BP}. Now, the hypothesis $\rm (H_1)$ is immediate because the diffusion process is derived from the local regular Dirichlet form $(\E\ui{\infty}, \D)$. Moreover, due to Theorem~\ref{BMI.thm20}, we confirm $\rm (H_2)$ and $\rm (H_3)$. Thus, we obtain the desired result by \cite[Theorem~8.1]{BP}.
\enddemo

\setcounter{equation}{0}

\section{Trace of the standard resistance form on $[0, \infty)$}\label{TRI}

In this section, we study the trace of the standard resistance form $(\E\ui{\infty}, \F\ui{\infty})$ on $I\ui{\infty} \times \{0\}$, where $I\ui{\infty} = [0, \infty)$, and the Dirichlet form associated with the trace. In particular, we are going to give an exact expression of the jump kernel of the trace and obtain estimates of the associated transition density in Theorem~\ref{TRI.thm10}, which includes Theorem~\ref{INT.thm10}. For simplicity, for a subset $A \subseteq [0, \infty)$, we identify $A$ with $A \times \{0\} \subseteq \BbR^2$. For example, we define $I\ui{m} = [0, 2^m] \subseteq [0, \infty)$ and use $I\ui{m}$ to denote $[0, 2^m] \times \{0\} \subseteq \BbR^2$ for $m \ge 0$.

\definition\label{TRI.def05}
Define $(\E\ui{\infty}_I, \F\ui{\infty}_I) = ([\E\ui{\infty}]_{I\ui{\infty}}, [\F\ui{\infty}]_{I\ui{\infty}})$. For $u \in \F\ui{\infty}$, define
\[
\E\ui{\infty}_{I, 0}(u) = \E\ui{\infty}_I(u) + u(0)^2.
\]
\enddefinition

\remark
The resistance metric associated with the trace $(\E\ui{\infty}_I, \F\ui{\infty}_I)$ is the restriction of $R\ui{\infty}$, which is the resistance metric associated with $(\E\ui{\infty}, \F\ui{\infty})$. For simplicity, we use $R\ui{\infty}$ to denote the resistance metric associated with $(\E\ui{\infty}_I, \F\ui{\infty}_I)$ as well.
\endremark

\thm\label{TRI.thm00}
$\F\ui{\infty}_I \cap C_0(I\ui{\infty})$ is dense subset of $\F\ui{\infty}_I$ with respect to the norm $\sqrt{\E\ui{\infty}_{I, 0}}$. In particular, 
\begin{equation}\label{TRI.eq10}
R\ui{\infty}(x, y) = \big(\inf\{\E\ui{\infty}_I(u)| u \in \F\ui{\infty} \cap C_0(I\ui{\infty}), u(x) = 1, u(y) = 0\}\big)^{-1}
\end{equation}
\endthm

\demo
Let $v \in \F\ui{\infty}_I$. Choose $u \in \F\ui{\infty}$ such that $u|_{I\ui{\infty}} = v$. Fix $\e > 0$. Then by Theorem~\ref{ISG.thm40} and its proof, there exists $u_* \in \F\ui{\infty} \cap C_0(K\ui{\infty})$ such that $\E\ui{\infty}(u - u_*) < \e$ and $u(p_0) = u_*(p_0)$. Let $v_* = u_*|_{I\ui{\infty}}$. Then
\[
\E\ui{\infty}_{I, 0}(v - v_*) = \E\ui{\infty}(h_{I\ui{\infty}}^{\E\ui{\infty}}(v - v_*)) \le \E\ui{\infty}(u - u_*) < \e.
\]
Thus we have shown that $\F\ui{\infty}_I \cap C_0(I\ui{\infty})$ is dense subset of $\F\ui{\infty}_I$. The equality \eqref{TRI.eq10} is straightforward.
\enddemo

Next we study the regular Dirichlet form derived from $(\E\ui{\infty}_I, \F\ui{\infty}_I)$.

\definition\label{ITE.def20}
(1)\,\,Define $J\ui{\infty}:(I\ui{\infty})^2\backslash \{(x, x)| x \in I\ui{\infty}\} \to [0, \infty)$ by
\begin{equation}\label{ITE.eq10}
J\ui{\infty}(x, y) = \frac{35}{16}\cdot\frac{14}{17}\Big(\frac {20}3\Big)^{n - 1}
\end{equation}
if $(x, y) \in \bigcup_{i \in \BbN}S_{n, i}$.\\
(2)\,\, Let $\nu$ be the Lebesgue measure on $I\ui{\infty}$. For $u \in \F\ui{\infty}_I \cap L^2(I\ui{\infty}, \nu)$, define
\[
\E_*(u) = \E\ui{\infty}_I(u) + \int_{I\ui{\infty}} |u(x)|^2\nu(dx).
\]
Moreover, let $\D\ui{\infty}_I$ be the closure of $\F\ui{\infty} \cap C_0(I\ui{\infty})$ with respect to the norm $\sqrt{\E_*(\cdot)}$.
\enddefinition

\remark
In Theorem~\ref{INT.thm10}, the definition \eqref{ITE.eq10} is stated as a claim of the theorem. Combining Definition~\ref{ITE.def20} and Theorem~\ref{TRI.thm10}, we see that there is no discrepancy between the two presentations of \eqref{ITE.eq10} after all.
\endremark

\remark
To make a comparison with $J\ui{m}$ defined later easier, we leave $\frac{35}{16}\cdot\frac{14}{17}$ in \eqref{ITE.eq10} not reduced. See Definition~\ref{TRI.def10}.
\endremark 

Now we have our main results on the Dirichlet forms associated with the trace.

\thm\label{TRI.thm10}
$(\E\ui{\infty}_I, \D\ui{\infty}_I)$ is a regular Dirichlet form on $L^2(I\ui{\infty}, \nu)$ and 
\begin{equation}\label{TRI.eq15}
\E\ui{\infty}_I(u) = \int_{I\ui{\infty}} (u(x) - u(y))^2J\ui{\infty}(x, y)\nu(dx)\nu(dy).
\end{equation}
for any $u \in \D\ui{\infty}_I$. Moreover, there exist an associated jointly continuous transition density $p\ui{\infty}_I: (0, \infty) \times I\ui{\infty} \times I\ui{\infty} \to [0, \infty)$, positive constants $c_3, c_4$ and $c_5$ such that
\begin{equation}\label{TRI.eq20}
p\ui{\infty}_I(t, x, y) \le c_3\min\bigg\{t^{-1/(\a + 1)}, \frac{t}{|x - y|^{\a + 2}}\bigg\}
\end{equation}
for any $(t, x, y) \in (0, \infty) \times I\ui{\infty} \times I\ui{\infty}$, and
\begin{equation}\label{TRI.eq30}
p\ui{\infty}_I(t, x, y) \ge c_4t^{-1/(\a + 1)}
\end{equation}
if $|x - y| \le c_5t^{\frac 1{1 + \a}}$.
\endthm

The estimates \eqref{TRI.eq20} and \eqref{TRI.eq30} are called the upper transition density estimate and the lower near diagonal transition density estimate, respectively. In \cite{ChenKi}, the following better off-diagonal lower estimate of $p\ui{\infty}(t, x, y)$ is obtained;
\[
p\ui{\infty}_I(t, x, y) \ge \begin{cases}
c_4t^{-\frac 1{1 + \a}}\quad&\text{if $|x - y| \le c_5t^{\frac 1{1 + \a}}$},\\
c_6J\ui{\infty}(x, y)t\quad&\text{otherwise},
\end{cases}
\]
See \cite{ChenKi} for details.\par
The rest of this section is devoted to a proof of Theorem~\ref{TRI.thm10}. To begin with, we identify the jump kernel $J\ui{m}(x, y)$ of the trace of $(\E\ui{m}, \F\ui{m})$ on $I\ui{m}$, which is one of the three boundary segments of the enlarged Sierpinski gasket $K\ui{m}$. Note that
\[
\sd{(I\ui{m})^2}{\{(x, x)| x \in I\ui{m}\}} = \bigsqcup_{n \ge 1 - m, 1 \le i \le 2^{n + m - 1}} S_{n, i}.
\]

\definition\label{TRI.def10}
For $m \ge 0$, define $J\ui{m}:\sd{(I\ui{m})^2}{\{(x, x)| x \in I\ui{m}\}} \to [0, \infty)$ by
 \[
 J\ui{m}(x, y) = \frac{35}{16}\Big(\frac{14}{17}\Big(\frac {20}3\Big)^{n - 1} + \frac{3}{17}\Big(\frac 3{20}\Big)^m\Big).
 \]
 if $(x, y) \in \bigcup_{1 \le i \le 2^{n + m - 1}} S_{n, i}$.
 \enddefinition
 
\lemma\label{TRI.lemma10}
Let $(\E^{(m)}_{I}, \F^{(m)}_{I})$ be the trace of $(\E^{(m)}, \F^{(m)})$ on $I\ui{m}$. Then
\begin{multline*}
\F^{(m)}_I \\= \Big\{u\Big| u \in C(I\ui{m}, d), \int_{I\ui{m} \times I\ui{m}} (u(x) - u(y))^2J\ui{m}(x, y)\nu(dx)\nu(dy) < \infty\Big\}
\end{multline*}
and
\[
\E^{(m)}_I(u) = \int_{I\ui{m} \times I\ui{m}} (u(x) - u(y))^2J\ui{m}(x, y)\nu(dx)\nu(dy)
\]
for any $u \in \F^{(m)}_I$.
\endlemma

\demo
In the case $m = 0$, this lemma is a part of \cite[Corollary~6.2]{KiKTaka}, where $J\ui{0}$ is denoted by $J_*$.
By the definition of $J\ui{m}$, 
\[
J\ui{m}(x, y) = \Big(\frac 3{20}\Big)^mJ\ui{0}(2^{-m}x, 2^{-m}y).
\]
By a scaling argument, we see that
\begin{align*}
\E^{(m)}_I(u, u) &= \Big(\frac 35\Big)^m\int_{I^2} (u(2^mx) - u(2^my))^2J\ui{0}(x, y)\nu(dx)\nu(dy)\\
&= \Big(\frac 3{20}\Big)^m\int_{(I\ui{m})^2} (u(X) - u(Y))^2J\ui{0}(2^{-m}X, 2^{-m}Y)\nu(dX)\nu(dY)\\
&= \int_{(I\ui{m})^2} (u(X) - u(Y))^2J\ui{m}(X, Y)\nu(dX)\nu(dY).
\end{align*}
\enddemo

\demo[Proof of Theorem~\ref{TRI.thm10}]
First we show that $(\E\ui{\infty}_I, \D\ui{\infty}_I)$ is a regular Dirichlet form on $L^2(I\ui{\infty}, \nu)$. Since $(\E\ui{\infty}, \F\ui{\infty})$ is regular and $I\ui{\infty}$ is closed, using \cite[Theorem~8.4]{Ki16}, we see that $(\E\ui{\infty}_I, \F\ui{\infty}_I)$ is regular. Then by \cite[Theorem~9.4]{Ki16}, if follows that $(\E\ui{\infty}_I, \D\ui{\infty}_I)$ is a regular Dirichlet form on $L^2(I\ui{\infty}, \nu)$. \par
Next we establish \eqref{TRI.eq15}. Let $v \in \F\ui{\infty}_I \cap C_0(I\ui{\infty})$. Choose $n \in \BbN$ such that $\supp{v} \subseteq I\ui{n}$. Let $m \ge n$. Set  $u_0 = h_{I\ui{\infty}}^{\E\ui{\infty}}(v)$ and $u_1 = h_{I\ui{m}}^{\E\ui{m}}(v|_{I\ui{m}})$. Then
\[
\E\ui{\infty}_I(v) = \E\ui{\infty}(u_0)\quad\text{and}\quad \E\ui{m}_I(v|_{I\ui{m}}) = \E\ui{m}(u_1)
\]
Since $(u_0|_{K\ui{m}})|_{I\ui{m}} = v|_{I\ui{m}}$, we see that
\[
\E\ui{m}(u_1) \le \E\ui{m}(u_0|_{I\ui{m}}) \le \E\ui{\infty}(u_0).
\]
By Lemma~\ref{TRI.lemma10}, we see that
\begin{equation}\label{TRI.eq100}
\int_{(I\ui{m})^2}(v(x) - v(y))^2J\ui{m}(x, y)\nu(dx)\nu(dy) = \E\ui{m}(u_1) \le \E\ui{\infty}_I(v).
\end{equation}
On the other hand, let $h_0 = h^{\E}_{V_0}(\chi_{p_0}|_{V_0})$, where $\chi_{p_0}$ is the characteristic function of the set $\{p_0\}$. Then $h_0:K \to [0, 1]$ and $\E(h_0) = 2$. Define $h_m:K\ui{\infty, m} \to [0, 1]$ by
\[
h_m(z) = \begin{cases}
h_0((\rho_2\ui{m})^{-1}(z))\quad&\text{if $z \in \rho_2\ui{m}(K)$,}\\
0\quad&\text{otherwise.}
\end{cases}
\]
Then $h_m \in \F\ui{\infty, m}$, $\displaystyle\E\ui{\infty, m}(h_m) = 2\Big(\frac 35\Big)^m$, and $h_m|_{[2^m, \infty) \times 0} = 0$.
Furthermore, define $u_2: K\ui{\infty} \to \BbR$ by
\[
u_2(z) = \begin{cases}
u_1(z)\quad&\text{if $z \in K\ui{m}$,}\\
u_1(p_2)^2h_m(z)\quad&\text{if $z \in K\ui{\infty, m}$}.
\end{cases}
\]
Since $u_1(p_1\ui{m}) = v(p_1\ui{m}) = 0$, we see that $u_2 \in \F\ui{\infty}$ and 
\[
\E\ui{\infty}(u_2) = \E\ui{m} + 2\Big(\frac 35\Big)^mu_1(p_2^{(m)})^2.
\]
Since $u_2|_{I\ui{\infty}} = v$, we see that $\E\ui{\infty}_I(v) \le \E\ui{\infty}(u_2)$. Therefore, the above equality yields
\begin{multline}\label{TRI.eq200}
\E\ui{\infty}_I(v) \\
\le \int_{(I\ui{m})^2}(v(x) - v(y))^2J\ui{m}(x, y)\nu(dx)\nu(dy) + 2\Big(\frac 35\Big)^mu_1(p_2^{(m)})^2.
\end{multline}
Now we consider what happens if we let $m \to \infty$. First, let $\overline{v}_m = \max_{x \in I\ui{m}} v(x)$ and $\underline{v}_m = \min_{x \in I\ui{m}} v(x)$ and define $\bar{u}_1: K\ui{m} \to \BbR$ by
\[
\bar{u}_1(z) = \begin{cases}
\overline{v}_m &\quad\text{if $u_1(z) \ge \overline{v}_m$,}\\
\underline{v}_m &\quad\text{if $u_1(z) \le \underline{v}_m$,}\\
u_1(z) &\quad\text{otherwise.}
\end{cases}
\]
Then the Markov property implies that $\bar{u}_1 \in \F\ui{m}$ and $\E\ui{m}(\bar{u}_1) \le \E\ui{m}(u_1)$. Since $\bar{u}_1|_{I\ui{m}} = v|_{I\ui{m}}$, Proposition~\ref{BRF.prop30} shows that $\bar{u}_1 = u_1$. Hence
\[
u_1(p_2\ui{m})^2 = |u_1(p_1\ui{m}) - u_1(p_2\ui{m})| \le \max_{x \in I\ui{\infty}} v(x) - \min_{x \in I\ui{\infty}} v(x).
\]
Therefore, as $m \to \infty$, 
\[
2\Big(\frac 35\Big)^mu_1(p_2^{(m)})^2 \to 0.
\]
Second, 
\begin{multline*}
 \int_{(I\ui{m})^2}(v(x) - v(y))^2J\ui{m}(x, y)\nu(dx)\nu(dy) \\=  \int_{(I\ui{m})^2}(v(x) - v(y))^2J\ui{\infty}(x, y)\nu(dx)\nu(dy)
 \\ + \frac{45}{272}\Big(\frac 3{20}\Big)^m\int_{(I\ui{m})^2}(v(x) - v(y))^2\nu(dx)\nu(dy).
 \end{multline*}
 The monotone convergence theorem implies that
 \begin{multline*}
 \int_{(I\ui{m})^2}(v(x) - v(y))^2J\ui{\infty}(x, y)\nu(dx)\nu(dy)\\ \to \int_{(I\ui{\infty})^2}(v(x) - v(y))^2J\ui{\infty}(x, y)\nu(dx)\nu(dy)
 \end{multline*}
 as $m \to \infty$. Moreover, as $m \to \infty$, 
 \begin{align*}
 &\Big(\frac 3{20}\Big)^m\int_{(I\ui{m})^2}(v(x) - v(y))^2\nu(dx)\nu(dy)\\
 = &\Big(\frac 3{20}\Big)^m\Bigg(2^{m + 1}\int_{I\ui{m}} v(x)^2\nu(dx) - 2\bigg(\int_{I\ui{m}} v(x)\nu(dx)\bigg)^2\Bigg) \to 0\\
 \end{align*}
 Thus we have shown 
 \begin{multline*}
 \int_{(I\ui{m})^2}(v(x) - v(y))^2J\ui{m}(x, y)\nu(dx)\nu(dy) \\
 \to \int_{(I\ui{\infty})^2} (v(x) - v(y))^2J\ui{\infty}(x, y)\nu(dx)\nu(dy).
 \end{multline*}
 as $m \to \infty$. Therefore, combining those estimates with \eqref{TRI.eq100} and \eqref{TRI.eq200}, we have shown that
 \[
 \E^{\ui{\infty}}_I(v) = \int_{(I\ui{\infty})^2} (v(x) - v(y))^2J\ui{\infty}(x, y)\nu(dx)\nu(dy).
 \]
 Next let $\{u_n\}_{n \ge 1}$ be a $\E_*$-convergent sequence in $\F\ui{\infty} \cap C_0(I\ui{\infty})$. Then $\{u_n\}_{n \ge 1}$ is also a convergent sequence in $L^2(I\ui{\infty}, \nu)$. Let $u$ be its limit. Choosing a subsequence, we may assume that $\{u_n(x)\}_{n \ge 1}$ converges to $u(x)$ as $n \to \infty$ for a.e $x \in I\ui{\infty}$. Define 
 \[
 \nu^2_J(A) = \int_{A} J\ui{\infty}(x, y)\nu(dx)\nu(dy)
 \]
 for a Borel set $A \subseteq (I\ui{\infty})^2$. Define $\Phi:\F\ui{\infty} \cap C_0(X) \to L^2((I\ui{\infty})^2, \nu^2_J)$ by $\Phi(u)(x, y) = u(x) - u(y)$. Then $\{\Phi{u_n}\}_{n \ge 1}$ is a convergent sequence in $L^2((I\ui{\infty})^2, \nu^2_J)$. Let $F$ be its limit. Since $u_n(x) - u_n(y)$ is convergent to $u(x) - u(y)$ for a.e. $(x, y) \in (I\ui{\infty})^2$, it follows that $F(x, y) = u(x) - u(y)$. Therefore, 
 \begin{multline*}
\E\ui{\infty}_I(u) = \lim_{n \to \infty} \E\ui{\infty}_I(u_n)\\ = \lim_{n \to \infty} \int_{(I\ui{\infty})^2}(u_n(x) - u_n(y))^2J\ui{\infty}(x, y)\nu(dx)\nu(dy) 
\\= \int_{(I\ui{\infty})^2}(u(x) - u(y))^2J\ui{\infty}(x, y)\nu(dx)\nu(dy).
\end{multline*}
Thus we have verified \eqref{TRI.eq15}. \par
The existence of the jointly continuous transition density is due to \cite[Theorem~10.4]{Ki16}. In \cite[Theorem~6.17]{GriHuLau2}, the estimates \eqref{TRI.eq20} and \eqref{TRI.eq30} have been shown to be equivalent to the combination of conditions ''parabolicity'', (V), ($J_{\le}), and (R)$ given in \cite{GriHuLau2}. Since $\nu$ is the Lebesgue measure, it follows that
\[
\frac 12r \le \nu(B(x, r)) \le r
\]
for any $x \in I\ui{\infty}$ and $r > 0$, where $B(x, r) = \{y| y \in I\ui{\infty}, |x - y| < r\}$. This immediately implies the condition (V). Since there exists $c > 0$ such that 
\[
J\ui{\infty}(x, y) \le \frac c{|x - y|^{\a + 2}}
\]
for any $x, y \in (I\ui{\infty})^2$ with $x \neq y$, we have the condition ($J_{\le}$). The condition (R) follows from \eqref{TRI.eq10} and \eqref{ISG.eq100}. Finally, note that
\[
|u(x) - u(y)|^2 \le R\ui{\infty}(x, y)\E\ui{\infty}_I(u) \le c_2|x - y|^{\a}\E\ui{\infty}_I(u).
\]
By \cite[Proposition~6.5]{GriHuLau2}, this and the condition (R) imply the ``parabolicity''. Thus, we have verified the estimates \eqref{TRI.eq20} and \eqref{TRI.eq30}.
\enddemo

\appendix

\setcounter{equation}{0}
\section{Basics of resistance forms}\label{APD}

In this appendix, we briefly review the theory of resistance forms and introduce several lemmas, Lemmas~\ref{BRF.lemma10} and \ref{BRF.lemma20} on resistance forms. One can find more detailed treatments on resistance forms in \cite{AOF, Ki11, Ki16}. The following definition of resistance forms is due to \cite{Ki16}. 

\definition\label{BRF.def10}
Let $X$ be a set. A pair $(\E, \F)$ is called a resistance form on $X$ if the following conditions (RF1) through (RF5) are satisfied:\\
(RF1)\,\,$\F$ is a linear subspace of $\ell(X)$ containing constants and $\E$ is a non-negative symmetric quadratic form on $\F$. $\E(f, f) = 0$ if and only if $f$ is a constant function on $X$.\\
(RF2)\,\,Let $\sim$ be an equivalence relation on $\F$ defined by $f \sim g$ if and only if $f - g$ is a constant function on $X$. Then $(\F/\sim, \E)$ is a real Hilbert space.\\
(RF3)\,\,If $x \neq y \in X$, then there exists $f \in \F$ such that $f(x) \neq f(y)$.\\
(RF4)\,\,For any $x, y \in X$, 
\[
\sup_{f \in \F, \E(f, f) \neq 0} \frac{|f(x) - f(y)|^2}{\E(f, f)} < \infty.
\]
(RF5)\,\,For any $f \in \F$, $\overline{f} \in \F$ and $\E(\overline{f}, \overline{f}) \le \E(f, f)$, where $\overline{f}$ is given by
\[
\overline{f}(x) = \begin{cases}
1 &\text{if $f(x) \ge 1$,}\\
0 &\text{if $f(x) \le 0$,}\\
f(x) & \text{if $0 < f(x) < 1$.}
\end{cases}
\]
For a resistance form $(\E, \F)$ on $X$, we denote the supremum in (RF4) by $R_{\E}(x, y)$ and call it the resistance metric associated with $(\E, \F)$.
\enddefinition

The name ``resistance metric'' is justified by the following proposition. 

\prop\label{BRF.prop10}
Let $(\E, \F)$ be a resistance form on a set $X$. Then $R_{\E}(\cdot, \cdot)$ is a metric on $X$.
\endprop

This fact is originally shown in \cite[Theorem~3.2]{Ki6}, where a resistance form was called a finite resistance form.

Next we introduce the notion of the trace of a resistance form according to \cite{Ki16}. We have slightly changed, however, the original notations in \cite{Ki16} to clarify the difference between a restriction and a trace.

\prop\label{BRF.prop30}
Let $(\E, \F)$ be a resistance form on a set $X$ and let $Y$ be a non-empty subset of $X$. Define
\[
[\F]_Y = \big\{u|_Y \big|u \in \F\big\}
\]
and
\[
[\E]_Y(v) = \inf\{\E(u, u)| u \in \F, u|_Y = v\}
\]
for $v \in [\F]_Y$. Then $([\E]_Y, [\F]_Y)$ is a resistance form on $Y$ and $R_{\E}(x, y) = R_{[\E]_V}(x, y)$ for any $x, y \in Y$. Moreover, there exists a linear map $h_Y^{\E}: [\F]_Y \to \F$ such that $h_Y^{\E}(v)|_Y = v$ and $h_Y^{\E}(v)$ is the unique element which satisfies
\[
[\E]_Y(v) = \E(h_Y^{\E_i}(v)).
\]
\endprop

The resistance form $([\E]_Y, [\F]_Y)$ on $Y$ is called the trace of the resistance form $(\E, \F)$ on $Y$.\par
Rewriting the definition of of the resistance metric in (RF4), we see that
\begin{equation}\label{BRF.eq50}
R(x, y) = \big(\inf\{\E(u, u)| u \in \F, u(0) = x, u(1) = y\})^{-1}.
\end{equation}
Comparing this with the definition of the trace, we obtain the following corollary.

\cor\label{BRF.cor10}
Let $(\E, \F)$ be a resistance form on $X$ and let $R(\cdot, \cdot)$ be the associated resistance metric. Then
\[
R(x, y) = \big([\E]_{\{x, y\}}(h_{\{x, y\}}^{\E}(\chi_x)\big)^{-1}
\]
for any $x, y \in X$ with $x \neq y$, where $\chi_x: X \to \BbR$ is given by $\chi_x(x) = 1$ and $\chi_x(y) = 0$ if $x \neq y$.
\endcor

Resistance forms on a finite set are identified with a weighted graph on the set as follows.

\prop\label{BRF.prop40}
Let $V$ be a finite set. Let $\F$ be a linear subspace of $\ell(V)$ and let $\E$ be a non-negative quadratic form on $\F$. The pair $(\E, \F)$ is a resistance form on $V$ if and only if $\F = \ell(V)$ and there exists $C: V \times V \to [0, \infty)$ such that $C(x, y) = C(y, x)$ for any $x, y \in V$, $(V, E_C)$ is a non-directed connected graph, where $E_C = \{(x, y)| x, y \in V, C(x, y) > 0\}$, and 
\[
\E(u, v)  = \frac 12\sum_{(x, y) \in E_C} C(x, y)(u(x) - u(y))(v(x) - v(y))
\]
for any $u, v \in \F$.
\endprop

The coefficient $C(x, y)$ appearing in the above proposition and its reciprocal correspond to the electrical conductance and the electrical resistance between $x$ and $y$, respectively. On the analogy of electrical circuits, $R_{\E}(x, y)$ is the effective resistance between $x$ and $y$.\par
In light of the above proposition, if $(\E, \F)$ is a resistance form on a finite set $V$, we say that $\E$ is a resistance form on $V$ or $(V, \E)$ is a resistance form.\par
The following two lemmas are useful for calculating a trace of a resistance form on a finite set.

\lemma\label{BRF.lemma10}
Let $V$ be a finite set and let $\E$ be a resistance form on $V$. Suppose that $V = \bigcup_{i = 1}^k V_k$, where $V_k$ is a non-empty subset of $V$ for any $k = 1, \ldots, k$, and that
\[
\E(u) = \sum_{i = 1}^k \E_i(u|_{V_i})
\]
for any $u \in \ell(V)$, where $\E_i$ is a resistance form on $V_i$ for any $i = 1, \ldots, k$. If 
\begin{equation}\label{BRF.eq100}
\bigcup_{i, j \in \{1, \ldots, k\}, i > j} (V_i \cap V_j) \subseteq U,
\end{equation}
then 
\[
[\E]_U(v) = \sum_{i = 1}^k [\E_i]_{U \cap V_i}(v|_{U \cap V_i})
\]
for any $v \in \ell(U)$.
\endlemma

\remark
 If $V_i$ is a single point, $\E_i$ is a trivial resistance form. We exclude such a case, i.e. we assume that $\#(V_i) \ge 2$ for any $i \in \{1, \ldots, k\}$.
 \endremark
 
 \demo
 For simplicity, we denote $U \cap V_i$, $h_U^{\E}$ and $h_{U \cap V_i}^{\E_i}$ by $U_i$, $h$ and $h_i$,  respectively.
 For $v \in \ell(U)$, note that $h(v)|_{U_i} = v|_{U_i}$, we see that
 \[
 [\E]_U(v) = \E(h(u)) = \sum_{i = 1}^k \E_i(h(u)|_{V_i}) \ge \sum_{i = 1}^k \E_i(h_i(v|_{U_i})) = \sum_{i = 1}^k [\E_i]_{U_i}(v|_{V_i})
 \]
 On the other hand, by \eqref{BRF.eq100}, if follows that
 \[
 h_i(v|_{U_i})|_{V_i \cap V_j} = v|_{V_i \cap V_j} = h_j(v|_{U_j})|_{V_i \cap V_j}.
 \]
 Hence there exists $u_* \in \ell(V)$ such that $u_*|_U = v$ and $u_*|_{V_i} = h_{U \cap V_i}^{\E_i}(v_i)$ for any $i \ge 1$. This yields
 \[
 \sum_{i = 1}^k [\E_i]_{U_i}(v|_{V_i}) = \E(u_*) \ge [\E]_U(v).
 \]
 Thus, we have shown the desired claim.
\enddemo

\lemma\label{BRF.lemma20}
{\rm (1) {\bf Series of resistors}:}\,\,\,Let $V = \{p_0, p_1, \cdots, p_k\}$, where $p_i \neq p_j$ if $i \neq j$ and let $R_1, \cdots R_k > 0$. Define a resistance form $\E$ on $V$ by
\[
\E(u) = \sum_{i = 1}^k \frac 1{R_i}(u(p_{i - 1}) - u(p_i))^2
\]
for $u \in \ell(V)$. Set $U = \{p_0, p_k\}$ and $R = \sum_{i = 1}^k R_i$. Then
\[
[\E]_{U}(v) = \frac 1R(v(p_0) - v(p_k))^2
\]
for any $v \in \ell(V_0)$.\\
{\rm (2) {\bf Parallel resistors}:}\,\, Let $R_1, \cdots R_k > 0$ and let $V = \{p, q\}$, where $p \neq q$. Define a resistance form $\E_i$ on $V$ by
\[
\E_i(u) = \frac 1{R_i}(u(p) - u(q))^2
\]
for each $i = 1, \ldots, k$. Set $R = \Big(\sum_{i = 1}^k \frac 1{R_i}\Big)^{-1}$. Then $\E = \sum_{i = 1}^k \E_k$ is a resistance form on $V$ and
\[
\E(u) = \frac 1R(u(p) - u(q))^2.
\]
\endlemma

Finally we introduce the notion of a compatible sequence of resistance forms. The following definitions and the theorem are originally given in \cite{Ki6}. More systematic treatment can be found in \cite{AOF}.

\definition\label{BRF.def20}
(1)\,\,Let $\{(V_m, \E_m)\}_{m \ge 0}$ be a sequence of pairs of a finite set $V_m$ and a resistance form $\E_m$ on $V_m$. The sequence $\{(V_m, \E_m)\}_{m \ge 0}$ is called a compatible sequence of resistance forms on finite sets if $V_m \subseteq V_{m + 1}$ and $[\E_{m + 1}]_{V_m} = \E_m$ for any $m \ge 0$.\\
(2)\,\,Let $\L = \{(V_m, \E_m)\}_{m \ge 0}$ be a compatible sequence of resistance forms on finite sets. Define
\[
V_{\L} = \bigcup_{m \ge 0} V_m,
\]
\[
\F_{\L} = \{u| u \in \ell(V_{\L}), \lim_{m \ge 0} \E_m(u|_{V_m}) < +\infty\}
\]
and
\[
\E_{\L}(u) = \lim_{m \ge \infty} \E_m(u|_{V_m})
\]
for $u \in \F_{\L}$. The pair $(\E_{\L}, \F_{\L})$ is called the limit of $\L$.
\enddefinition

\thm\label{BRF.thm10}
Let $\L = \{(V_m, \E_m)\}_{m \ge 0}$ be a compatible sequence of resistance forms on finite sets. Then $(\E_L, \F_L)$ is a resistance form on $V_{\L}$. Let $R_{\L}$ be the resistance metric associated with $(\E_{\L}, \F_{\L})$. Then 
\[
R_{\L}(x, y) = R_{\E_m}(x, y)
\]
for any $x, y \in V_m$. Moreover let $(X, R)$ be the completion of the metric space $(V_{\L}, R_{\L})$. Then any $u \in \F_{\L}$ is extended to a continuous function on $X$ and $\F_{\L}$ is identified with a subset of continuous functions on $(X, R)$. Through this extension, $(\E_{\L}, \F_{\L})$ is a resistance form on $X$ and $R$ is the associated resistance metric.
\endthm

\end{document}